\documentclass{svmult}
\usepackage{amsmath,amsfonts,amssymb,verbatim}

\spnewtheorem{theorem*}{Theorem}[section]{\bf}{\it} 
\spnewtheorem{lemma*}[theorem*]{Lemma}{\bf}{\it}
\spnewtheorem{proposition*}[theorem*]{Proposition}{\bf}{\it}
\spnewtheorem{corollary*}[theorem*]{Corollary}{\bf}{\it}
\spnewtheorem{definition*}[theorem*]{Definition}{\bf}{\it}
\spnewtheorem{example*}[theorem*]{Example}{\bf}{\it}
\spnewtheorem{remark*}[theorem*]{Remark}{\bf}{\it}

\def\bi{\tilde\imath}
\def\bj{\tilde\jmath}

\def\CC{{\mathbb C}}
\def\com{\ts,\hskip-.5pt}

\def\d{\partial}
\def\de{\delta}
\def\De{\Delta}

\def\End{\operatorname{End}\ts}
\def\ep{\varepsilon}

\def\f{\mathfrak{f}}

\def\g{\mathfrak{g}}
\def\ge{\geqslant}
\def\gl{\mathfrak{gl}}

\def\h{\mathfrak{h}}

\def\ka{\kappa}

\def\La{{\rm\Lambda}}
\def\la{\lambda}
\def\lac{\check\la}
\def\lcd{\ts,\ldots,}
\def\le{\leqslant}

\def\muc{\check\mu}

\def\ns{\hskip-1pt}

\def\ot{\otimes}

\def\Pb{\,\overline{\hskip-2.5pt P\hskip-4.5pt\phantom{\bar{t}}}\ts}
\def\Ph{\Phi}
\def\ph{\varphi}
\def\Phd{\dot\Phi}
\def\Pt{\widetilde{P}}

\def\Rb{\,\overline{\hskip-2.5pt R\hskip-4.5pt\phantom{\bar{t}}}\ts}
\def\Rh{\widehat{R}}
\def\rhoc{\check\rho}
\def\Rt{\widetilde{R}}

\def\sgn{\operatorname{sign}}
\def\so{\mathfrak{so}}
\def\sp{\mathfrak{sp}}
\def\SY{\operatorname{SY}}

\def\t{\mathfrak{t}}
\def\th{\theta}
\def\ts{\hskip1pt}
\def\Tt{\widetilde{T}}
\def\tts{\hskip.5pt}

\def\U{\operatorname{U}}

\def\Y{\operatorname{Y}}

\def\ZZ{{\mathbb Z}}

\begin{document}

\title*{Irreducible representations of Yangians}
\titlerunning{Irreducible representations of Yangians}

\author{Sergey Khoroshkin\inst{1,2},
Maxim Nazarov\inst{3} and Paolo Papi\inst{4}} 
\authorrunning{Khoroshkin, Nazarov and Papi}

\institute{\!Institute for Theoretical and
Experimental Physics, Moscow 117259, Russia 
\and\!Department of Mathematics,
Higher School of Economics,
Moscow 117312, Russia
\and 
\!Department of Mathematics, 
University of York, 
York YO10 5DD, England
\and
\!Dipartimento di Matematica,
Sapienza Universit\`a di Roma,
Roma 00185, Italy}

 
\maketitle
\smartqed
\renewcommand{\theequation}{\thesection.\arabic{equation}}
\makeatletter 
\@addtoreset{equation}{section} 
\makeatother


\section*{\normalsize\bf 0. Introduction}
\label{S0}
\setcounter{section}{0}
\setcounter{theorem*}{0}
\setcounter{equation}{0}

Take the general linear Lie algebra $\gl_n$ over the complex field $\CC$
and consider the corresponding polynomial current Lie algebra $\gl_n[x]\ts$.
The Yangian $\Y(\gl_n)$ is a deformation of the
universal enveloping algebra $\U(\gl_n[x])$
in the class of Hopf algebras. Irreducible finite-dimensional
representations of~the associative algebra $\Y(\gl_n)$ were
classified by Drinfeld \cite{D2}. In Subsections 1.1 and 1.2
of the present article we recall this classification, together with
the original definition of the Yangian $\Y(\gl_n)$ coming from the theory 
of quantum integrable systems.

The additive group $\CC$ acts on $\Y(\gl_n)$ by Hopf algebra
automorphisms, see \eqref{tauz} for the definition of automorphism
corresponding to the element $t\in\CC\ts$. The associative algebra 
$\Y(\gl_n)$ contains $\U(\gl_n)$ as a subalgebra, and
admits a homomorphism onto 
$\U(\gl_n)$ identical on this subalgebra\tts; see \eqref{eval}.
By pulling the representation of $\U(\gl_n)$ on an exterior power of
$\CC^n$ back through the homomorphism \eqref{eval}, and then
back through any automorphism \eqref{tauz}, we get
an irreducible representation of $\Y(\gl_n)$ called \textit{fundamental\/}.

The work \cite{D2} provided only a parametrization
of all irreducible finite-dimensional
representations of $\Y(\gl_n)\ts$. The results of 
Akasaka and Kashiwara \cite{AK} showed that 
up to an automorphism of $\Y(\gl_n)$ of the form
\eqref{fut}, each of these representations
arises as a quotient of a tensor product of fundamental ones\tts;
see also the earlier works
of Chari and Pressley \cite{CP} and Cherednik \cite{C2}.
This fact can also be derived from the results of
\text{Nazarov and Tarasov~\cite{NT1}\tts;}
further results
were obtained by Chari \cite{C}.
Note that the works \cite{AK} and \cite{C} dealt with
representations of quantum affine algebras. For
connections to the representation theory of Yangians
see the work of Molev, Tolstoy and Zhang \cite{MTZ}
and the more recent work of Gautam and Toledano-Laredo~\cite{GT}.

Now let $\g_n$ be one of the two Lie algebras $\so_n\com\sp_n\ts$.
Regard $\g_n$ as the Lie subalgebra of $\gl_n$ 
preserving a non-degenerate bilinear form on 
the vector space $\CC^n$,
symmetric in the case $\g_n=\so_n\ts$, 
or alternating in the case $\g_n=\sp_n\ts$. 
For any $n\times n$ matrix $X$ let $\widetilde{X}$ 
be the conjugate of $X$ relative to this form.
As an associative algebra, the twisted Yangian 
$\Y(\g_n)$ is a deformation of the universal
enveloping algebra of the \textit{twisted\/} polynomial current Lie algebra
$$
\{X(x)\in\gl_n[x]:X(-x)=-\ts\widetilde{X}(x)\}\,.
$$
This is not a Hopf algebra deformation. However,
$\Y(\g_n)$ is a one-sided coideal
subalgebra in the Hopf algebra $\Y(\gl_n)\ts$.
Moreover, 
$\Y(\g_n)$ contains $\U(\g_n)$ as a subalgebra, and
has a homomorphism onto 
$\U(\g_n)$ identical on this subalgebra. 

The definition of the twisted Yangian $\Y(\g_n)$
was given by Olshanski~\cite{O} with some help
from the second author of the present article.
This definition was motivated by the works of Cherednik \cite{C1} and
Sklyanin \cite{S} on the quantum integrable systems with boundary conditions.
Irreducible finite-dimensional
representations of the associative algebra $\Y(\g_n)$ have been
classified by Molev, see the recent book \cite{M}. 
In Subsections 1.4 and 1.5 we recall this classification, 
together with the definition of the twisted Yangian $\Y(\g_n)\ts$.

A new approach to the representation theory
of the Yangian $\Y(\gl_n)$ and of its twisted analogues
$\Y(\so_n)\com\Y(\sp_n)$ was developed by Khoroshkin, 
Nazarov and Vinberg in
\cite{KN,KNV}. In particular, it was proved in \cite{KN}
that up to an automorphism of $\Y(\sp_n)$ of the form \eqref{ffut},
any irreducible finite-dimensional representation of $\Y(\sp_n)$
is a quotient of some tensor product of fundamental
representations of $\Y(\gl_n)\ts$. 
The tensor product is regarded as
a representation of $\Y(\sp_n)$ by restriction from $\Y(\gl_n)\ts$.
Similar result was also proved in \cite{KN} for those representations
of $\Y(\so_n)$ which integrate from $\so_n\subset\Y(\so_n)$
to the special orthogonal group $\mathrm{SO}_n\ts$.
Moreover, the work \cite{KN} provided new proofs of the above mentioned
results of \cite{AK} for $\Y(\gl_n)\ts$.
We summarize the results of \cite{KN} in Subsections~1.3~and~1.6
of the present article.

The realizations of irreducible representations of $\Y(\gl_n)$ and $\Y(\g_n)$
in \cite{AK} and \cite{KN} were not quite explicit. 
To make them more explicit is
the main aim of the present article.
In Subsections 2.6 and 2.7 we give explicit formulas
for intertwining operators of tensor 
products of fundamental representations of $\Y(\gl_n)\ts$,
such that the quotients by the kernels of these operators realize
all irreducible finite-dimensional
representations of $\Y(\gl_n)\ts$,
up to automorphisms \eqref{fut}. 
In Subsections 3.4, 3.5, 3.6
we give analogues of these formulas for $\Y(\g_n)\ts$.  

To obtain the irreducible representations
of $\Y(\gl_n)$ and $\Y(\g_n)$ as quotients of tensor products of
fundamental representations, in \cite{KN} 
we imposed certain conditions on the tensor products\ts;
see Theorems \ref{1.1} and \ref{1.2} below.
Some of these conditions are not necessary.
For 
$\Y(\gl_n)$ the condition \eqref{chercon}
did not appear in \cite{AK} and \cite{NT1}, but did appear
in \cite{C2}. In Subsection~2.8 we remove this condition from
Theorem~\ref{1.1} by using the results of Subsection 2.5.
In Subsection~3.7
we remove the conditions
\eqref{con1},\eqref{con2},\eqref{con3} from Theorem~\ref{1.2} for 
$\Y(\g_n)$ by using the results of Subsection 3.3.
Thus we extend the results of \cite{KN} for both 
$\Y(\gl_n)$ and $\Y(\g_n)\ts$,
up to the level of \cite{AK} and \cite{NT1} for $\Y(\gl_n)\ts$.


\section*{\normalsize\bf 1. Representations of Yangians}
\label{S1}
\setcounter{section}{1}
\setcounter{theorem*}{0}
\setcounter{equation}{0}


\textbf{1.1.}
First consider the \textit{Yangian\/} $\Y(\gl_n)$
of the Lie algebra $\gl_n\ts$.
This is a complex unital associative algebra
with a family of generators
$
T_{ij}^{\ts(1)},T_{ij}^{\ts(2)},\ts\ldots
$
where $i\com j=1\lcd n\ts$.
Defining relations for them
can be written using the series
\begin{equation*}
T_{ij}(x)=
\de_{ij}+T_{ij}^{\ts(1)}x^{-\ns1}+T_{ij}^{\ts(2)}x^{-\ns2}+\,\ldots
\end{equation*}
where $x$ is a formal parameter. Let $y$ be another formal parameter.
Then the defining relations in the associative algebra $\Y(\gl_n)$
can be written as
\begin{equation}
\label{yrel}
(x-y)\,[\ts T_{ij}(x)\ts,T_{kl}(y)\ts]\ts=\;
T_{kj}(x)\ts T_{il}(y)-T_{kj}(y)\ts T_{il}(x)\,.
\end{equation}

The algebra $\Y(\gl_1)$ is commutative. 
By \eqref{yrel} for any $t\in\CC$
the assignments
\begin{equation}
\label{tauz}
T_{ij}(x)\ts\mapsto\,T_{ij}(x-t)
\end{equation}
define an automorphism of the algebra $\Y(\gl_n)\ts$.
Here each of the formal
power series $T_{ij}(x-t)$ in $(x-t)^{-1}$ should be re-expanded
in $x^{-1}$. Every assignment \eqref{tauz} is a correspondence
between the respective coefficients of series in $x^{-1}$.
Relations \eqref{yrel} also show that for any
formal power series $f(x)$ in $x^{-1}$ with coefficients from
$\CC$ and leading term $1$, the assignments
\begin{equation}
\label{fut}
T_{ij}(x)\ts\mapsto\,f(x)\,T_{ij}(x)
\end{equation}
define an automorphism of the algebra $\Y(\gl_n)\ts$.
The subalgebra consisting of all elements of $\Y(\gl_n)$ 
which are invariant under every automorphism 
of the form \eqref{fut}, is called the 
\textit{special Yangian\/} of $\gl_n\ts$, and is denoted by
$\SY(\gl_n)\ts$. Two representations of the 
algebra $\Y(\gl_n)$ are called $\textit{similar\/}$ if they differ
by an automorphism of the form \eqref{fut}.

Let $E_{ij}\in\gl_n$ be the standard matrix units.
By \eqref{yrel} the assignments
\begin{equation}
\label{eval}
T_{ij}(x)\ts\mapsto\,\de_{ij}+E_{ij}\ts x^{-1}
\end{equation}
define a homomorphism of unital associative algebras
$\Y(\gl_n)\to\U(\gl_n)\ts$.
There is also an embedding $\U(\gl_n)\to\Y(\gl_n)\ts$,
defined by mapping $E_{ij}\mapsto T_{ij}^{\ts(1)}$. So
$\Y(\gl_n)$ contains the universal enveloping
algebra $\U(\gl_n)$ as a subalgebra. The homomorphism \eqref{eval} is
evidently identical on the subalgebra $\U(\gl_n)\subset\Y(\gl_n)\ts$.

The Yangian $\Y(\gl_n)$ is a Hopf algebra over the field $\CC\ts$.
The comultiplication $\De:\Y(\gl_n)\to\Y(\gl_n)\ot\Y(\gl_n)$ is defined by
the assignments
\begin{equation}
\label{1.33}
\De:\,T_{ij}(x)\ts\mapsto\ts\sum_{k=1}^n\ T_{ik}(x)\ot T_{kj}(x)\,.
\end{equation}
When taking tensor products of $\Y(\gl_n)\ts$-modules,
we will use \eqref{1.33}. 
The counit homomorphism
$\Y(\gl_n)\to\CC$ is defined by the assignments
$
T_{ij}(x)\ts\mapsto\ts\de_{ij}\ts.
$
By using this homomorphism, one defines the 
\textit{trivial} representation of $\Y(\gl_n)\ts$.

Further,
let $T(x)$ be the $n\times n$ matrix
whose $i\com j$ entry is the series 
$T_{ij}(x)\ts$.
The antipodal map $\Y(\gl_n)\to\Y(\gl_n)$ is defined by mapping
$
T(x)\mapsto T(x)^{-1}.
$
Here each entry of the inverse matrix $T(x)^{-1}$
is a formal power series in $x^{-1}$ with coefficients
from the algebra $\Y(\gl_n)\ts$, and the assignment
$
T(x)\mapsto T(x)^{-1}
$
is a correspondence between the respective matrix entries.

The special Yangian $\SY(\gl_n)$ is a Hopf subalgebra
of $\Y(\gl_n)\ts$. Moreover, it is isomorphic to
the Yangian $\Y(\mathfrak{sl}_n)$ of the special linear Lie
algebra $\mathfrak{sl}_n\subset\gl_n\ts$ studied in
\cite{D1,D2}. For the proofs of these facts see
\cite[Subsection 1.8]{M}.


\bigskip\noindent
\textbf{1.2.}
Up to their equivalence and similarity, the irreducible
finite-dimensional representations of the associative algebra
$\Y(\gl_n)$ are parametrized by sequences of $n-1$
monic polynomials $P_1(x)\lcd P_{n-1}(x)$ with complex coefficients.
In particular, $P_1(x)=\ldots=P_{n-1}(x)=1$ for the trivial
representation of $\Y(\gl_n)\ts$.

This paramet\-rization was given by Drinfeld \cite[Theorem 2]{D2}.
In the present article we will use another parametrization. It 
can be obtained by combining the results of Arakawa and Suzuki \cite{AS}
with those of \cite{D1}. Namely, consider the Lie algebra $\gl_m\ts$. 
Let $\t_m$ be a Cartan subalgebra of $\gl_m\ts$.
Up to equivalence and similarity, we will parametrize the 
non-trivial irreducible
finite-dimensional representations of $\Y(\gl_n)$
by $m=1\com2\com\ts\ldots$ and
by certain orbits in $\t_m^*\times\t_m^*$
under the diagonal shifted action of the Weyl group of $\gl_m\ts$.

Pick the Cartan subalgebra $\t_m$
of $\gl_m$ with the basis $E_{11}\lcd E_{mm}$ of the diagonal matrix units. 
For any weight $\la\in\t_m^*$ define the sequence
$\la_1\lcd\la_m$ of its {labels} by setting
$\la_a=\la\ts(E_{aa})$ for $a=1\lcd m\ts$.
In particular, for the half-sum $\rho\in\t_m^*$
of the positive roots we have $\rho_a=m/2-a+1/2\ts$. 
The Weyl group of $\gl_m$ is the symmetric group $S_m\ts$.
It acts on the Cartan subalgebra $\t_m$ by permuting the basis vectors 
$E_{11}\lcd E_{mm}\ts$.
Hence it acts on any weight $\la\in\t_m^*$ by permuting its labels. 
The {shifted} action of $w\in S_m$ on $\t_m^*$ is given~by
\begin{equation}
\label{shiftact}
w\circ\la=w\ts(\la+\rho)-\rho\ts.
\end{equation}

For our parametrization we will be using only
the orbits of the pairs $(\la\com\mu)\in\t_m^*\times\t_m^*$
such that all labels of the weight $\nu=\la-\mu$ belong to
$\{\tts1\lcd n-1\}\ts$. 
Given any such a pair $(\la\com\mu)$ define
a sequence $P_1(x)\lcd P_{n-1}(x)$ of monic polynomials
with complex coefficients, as follows. For each $i=1\lcd n-1$ put
\begin{equation}
\label{pi}
P_i(x)=\prod_{\nu_a=i}\,(x-\mu_a-\rho_a)
\end{equation}
where $\mu_1\lcd\mu_m$ and $\nu_1\lcd\nu_m$ are
the labels of $\mu$ and $\nu$ respectively.
Note that the simultaneous shifted action of the group $S_m$ on
the weights $\la$ and $\mu$ gives the usual permutational action
of $S_m$ on the labels of $\mu+\rho$ and of $\nu\ts$. Therefore
each polynomial \eqref{pi} depends
only on the $S_m\ts$-orbit of the pair
$(\la\com\mu)\in\t_m^*\times\t_m^*\ts$.
Moreover, the orbit is determined by the polynomials 
\eqref{pi} uniquely.
Furthermore, any sequence of monic polynomials $P_1(x)\lcd P_{n-1}(x)$
of the total degree $m$ 
with complex coefficients arises in this way.

In the next subsection, to each of these orbits
we will attach an irreducible $\Y(\gl_n)\ts$-module.
Its Drinfeld polynomials
$P_1(x)\lcd P_{n-1}(x)$ will be given by \eqref{pi}.
For the definition of the Drinfeld polynomials
of an arbitrary irreducible finite-dimensional $\Y(\gl_n)\ts$-module
used here see \cite[Subsection 5.1]{KN}.

 
\bigskip\noindent
\textbf{1.3.}
For $k=0\com1\lcd n$ consider
the exterior power $\La^k(\CC^n)$
of the defining $\gl_n$-module $\CC^n\ts$.
Using the homomorphism \eqref{eval}, regard it as
a module over the Yangian $\Y(\gl_n)\ts$. For $t\in\CC$
denote by $\Ph^{\ts k}_{\tts t}$ the $\Y(\gl_n)\ts$-module obtained by
pulling the $\Y(\gl_n)\ts$-module $\La^k(\CC^n)$ back through the
automorphism \eqref{tauz}.

Now take any $(\la\com\mu)\in\t_m^*\times\t_m^*$
such that all labels of the weight $\nu=\la-\mu$ belong to
the set $\{\tts1\lcd n-1\}\ts$. Let the weights 
$\la$ and $\mu$ vary so that $\nu$ is fixed.
For the Cartan subalgebra $\t_m$ of $\gl_m$
the weight $\mu$ is \textit{generic\/} 
if $\mu_a-\mu_b\notin\ZZ$ whenever $a\neq b\ts$.
Consider the tensor product of $\Y(\gl_n)\ts$-modules
\begin{equation}
\label{basprod1}
\Ph^{\,\nu_1}_{\mu_1+\rho_1+\frac12}
\ot\ts\ldots\ts\ot\ts
\Ph^{\,\nu_m}_{\mu_m+\rho_m+\frac12}\,.
\end{equation}
It is known that the $\Y(\gl_n)\ts$-module \eqref{basprod1}
is irreducible if (but not~only~if) 
the weight $\mu$ is generic,
see \cite[Theorem 4.8]{NT2} for a more general result.
Moreover, if (but not only if) $\mu$ is generic then 
all the $\Y(\gl_n)\ts$-modules obtained from \eqref{basprod1}
by permuting the tensor factors, are 
equivalent to each other. In particular,
for every generic $\mu$
there is a unique non-zero $\Y(\gl_n)\ts$-intertwining operator
\begin{equation}
\label{intoper1}
\Ph^{\,\nu_1}_{\mu_1+\rho_1+\frac12}
\ot\ts\ldots\ts\ot\ts
\Ph^{\,\nu_m}_{\mu_m+\rho_m+\frac12}
\to\,
\Ph^{\,\nu_m}_{\mu_m+\rho_m+\frac12}
\ot\ts\ldots\ts\ot\ts
\Ph^{\,\nu_1}_{\mu_1+\rho_1+\frac12}
\end{equation}
corresponding to the permutation of tensor factors 
of the maximal length. We will
denote this operator by $I(\mu)\ts$;
it is unique up to a multiplier from $\CC\,\backslash\ts\{0\}\ts$.

For all generic weights $\mu\ts$,
the source vector spaces of the operators
$I(\mu)$ are the same. The target vector spaces
of all $I(\mu)$ also coincide with each other. 
Hence $I(\mu)$ is a function of
$\mu$ taking values in the space of linear operators
\begin{equation*}
\label{extoper1}
\La^{\nu_1}(\CC^n)
\ot\ts\ldots\ts\ot\ts
\La^{\nu_m}(\CC^n)
\,\to\,
\La^{\nu_m}(\CC^n)
\ot\ts\ldots\ts\ot\ts
\La^{\nu_1}(\CC^n)
\end{equation*}
between tensor products of exterior powers of $\CC^n\ts$.
The multipliers from $\CC\,\backslash\{0\}$
can be chosen so that $I(\mu)$ becomes a rational function of
$\mu\ts$, see Section 2 for a particular choice. Any such a choice
allows to determine intertwining operators \eqref{intoper1}
for those non-generic $\mu$ where 
$I(\mu)$ is regular.
Thus the weight $\mu$ need not be generic anymore, 
so that the $\Y(\gl_n)\ts$-module \eqref{basprod1} may be reducible.
This way of determining intertwining operators
between $\Y(\gl_n)\ts$-modules goes back to the work of Cherednik
\cite{C2} and is commonly called the~\text{\it fusion procedure}.

For the Cartan subalgebra $\t_m$ of $\gl_m$
the weight $\la$ will be called \textit{dominant} 
if $\la_a-\la_b\neq-\,1\com-\,2\com\ts\ldots$ whenever $a<b\ts$.
The pair $(\la\com\mu)\in\t_m^*\times\t_m^*$
will be called \textit{good\/} if the weight $\la+\rho$
is dominant and moreover
\begin{equation}
\label{chercon}
\nu_a\ge\nu_b
\quad\text{whenever}\quad
\la_a+\rho_a=\la_b+\rho_b
\quad\text{and}\quad
a<b\ts.
\end{equation}
The orbit of any $(\la\com\mu)\in\t_m^*\times\t_m^*$
under the shifted action of $S_m$ on $\t_m^*$
does contain a good pair. Here we assume only that
$\nu_1\lcd\nu_m\in\{\tts1\lcd n-1\}\ts$.
In its present form, the next theorem
was proved in \cite[Subsection~5.1]{KN}.
It will be generalized in Subsection 2.8 of the present article.

\begin{theorem*}
\label{1.1}
For the fixed weight $\nu=\la-\mu$
the multipliers from $\CC\,\backslash\{0\}$ of the 
operators \eqref{intoper1} for all generic weights $\mu$ can be chosen so 
that the rational function $I(\mu)$ is regular and 
non-zero whenever the pair $(\la\com\mu)$ is good.
Then for any good pair $(\la\com\mu)$ the quotient of \eqref{basprod1}
by the kernel of operator $I(\mu)$
is an irreducible $\Y(\gl_n)\ts$-module with the
Drinfeld polynomials given by \eqref{pi}.
\end{theorem*}

 
\smallskip\noindent
\textbf{1.4.}
Let $\g_n$ be one of the two Lie algebras 
$\so_n\com\sp_n\ts$.
We will regard $\g_n$ as the Lie subalgebra of $\gl_n$ 
preserving a non-degenerate bilinear form $\langle\ ,\,\rangle$ on $\CC^n$,
symmetric in the case $\g_n=\so_n\ts$, 
or alternating in the case $\g_n=\sp_n\ts$. 
In the latter case, the positive integer $n$ has to be even. 
When considering $\so_n\com\sp_n$ simultaneously
we will use the following convention. Whenever
the double sign $\ts\pm\ts$ or $\ts\mp\ts$ appears, 
the upper sign will correspond to the case
of a symmetric form on $\CC^{\ts n}$ so that $\g_n=\so_n\ts$.
The lower sign will correspond to the case
of an alternating form on $\CC^{\ts n}$ so that $\g_n=\sp_n\ts$.

Let $\Tt(x)$ be the conjugate of the matrix $T(x)$
relative to the form~$\langle\ ,\ts\rangle$~on~$\CC^{n}$.
An involutive automorphism of the algebra $\Y(\gl_n)$
is defined by the assignment
\begin{equation}
\label{transauto}
T(x)\mapsto\Tt(-x)\ts.
\end{equation}
This assignment is a correspondence between respective matrix entries.
Now consider the matrix product $S(x)=\Tt(-x)\,T(x)\ts$. Its $ij$ entry
is the series
\begin{equation}
\label{yser}
S_{ij}(x)=
\sum_{k=1}^n\,\Tt_{ik}(-\ts x)\,T_{kj}(x)=
\de_{ij}+S_{ij}^{\ts(1)}x^{-\ns1}+S_{ij}^{\ts(2)}x^{-\ns2}+\,\ldots
\end{equation}
with coefficients from the algebra $\Y(\gl_n)\ts$.
The \textit{twisted Yangian\/} corresponding to the Lie algebra $\g_n$
is the subalgebra of $\Y(\gl_n)$ generated by the coefficients
$
S_{ij}^{\ts(1)},S_{ij}^{\ts(2)},\ts\ldots
$
where $i\com j=1\lcd n\ts$.
We denote this subalgebra by $\Y(\g_n)\ts$.

The algebras $\Y(\so_n)$ corresponding to different choices
of the symmetric form $\langle\ ,\,\rangle$ on $\CC^n$
are isomorphic to each other, and so are
the algebras $\Y(\sp_n)$ corresponding to different choices
of the alternating form $\langle\ ,\,\rangle$ on $\CC^n\ts$.
These isomorphisms can be described explicitly, see
for instance \cite[Corollary 2.3.2]{M}.

The subalgebra $\Y(\g_n)\ts\cap\,\SY(\gl_n)$ of $\Y(\gl_n)$
is denoted by $\SY(\g_n)\ts$, and is called the
\textit{special twisted Yangian\/} corresponding to $\g_n$.
The automorphism \eqref{fut} of $\Y(\gl_n)$
determines an automorphism of $\Y(\g_n)$ which maps
\begin{equation}
\label{ffut}
S(x)\mapsto f(x)\,f(-x)\,S(x)\ts.
\end{equation}
The subalgebra $\SY(\g_n)$ of $\Y(\g_n)$ consists of
the elements fixed by all such automorphisms. 
Two representations of the algebra 
$\Y(\g_n)$ are called $\textit{similar\/}$ if they differ
by such an automorphism.

There is an analogue for $\Y(\g_n)$ of the
homomorphism $\Y(\gl_n)\to\U(\gl_n)$ defined by \eqref{eval}.
Namely, one can define a homomorphism $\Y(\g_n)\to\U(\g_n)$ by 
\begin{equation}
\label{pin}
S_{ij}(x)\,\mapsto\,\de_{ij}+
\frac{E_{ij}-\widetilde{E}_{ij}}
{\textstyle x\pm\frac12}
\end{equation}
where $\widetilde{E}_{ij}$ is the conjugate of
the matrix unit $E_{ij}\in\gl_n$
relative to the form $\langle\ ,\,\rangle$ on $\CC^{\ts n}\ts$.
This can be proved by using the defining relations for the generators 
$
S_{ij}^{\ts(1)},S_{ij}^{\ts(2)},\ts\ldots
$
which we do not reproduce here\ts;
see \cite[Proposition~2.1.2]{M} for the proof.
Further, there is an embedding $\U(\g_n)\to\Y(\g_n)$ defined by mapping
$$
E_{ij}-\widetilde{E}_{ij}\mapsto S_{ij}^{\ts(1)}\ts.
$$
Hence the twisted Yangian
$\Y(\g_n)$ contains the universal enveloping
algebra $\U(\g_n)$ as a subalgebra.
The homomorphism $\Y(\g_n)\to\U(\g_n)$ defined by \eqref{pin}
is evidently identical on the subalgebra $\U(\g_n)\subset\Y(\g_n)\ts$.

The twisted Yangian $\Y(\g_n)$ is not only a subalgebra of
$\Y(\gl_n)\ts$, it is also a right coideal of the coalgebra
$\Y(\gl_n)$ relative to the comultiplication \eqref{1.33}. Indeed,
by the definition of the series \eqref{yser} we get
$$
\De\ts(\,S_{ij}(x)))\,=
\sum_{k,l=1}^n\,S_{kl}(x)
\ts\ot\ts
\Tt_{ik}(-\ts x)\,T_{lj}(x)\ts.
$$
Therefore
$$
\De\ts(\ts\Y(\g_n))\subset\Y(\g_n)\ot\Y(\gl_n)\ts.
$$
Hence by taking a tensor product of an $\Y(\g_n)\ts$-module with
an $\Y(\gl_n)\ts$-module we get another $\Y(\g_n)\ts$-module.

The \textit{trivial\/} representation of the
algebra $\Y(\g_n)$ is defined by restricting the counit homomorphism 
$\Y(\gl_n)\to\CC$ to the subalgebra $\Y(\g_n)\subset\Y(\gl_n)\ts$.
Under this representation $S_{ij}(x)\mapsto\de_{ij}\ts$. Note that
restricting any representation of $\Y(\gl_n)$
to the subalgebra $\Y(\g_n)$ amounts to
taking the tensor product of that representation
of $\Y(\gl_n)$ with the
trivial representation of the algebra $\Y(\g_n)\ts$.

A parametrization of irreducible finite-dimensional
representations of the algebra $\Y(\g_n)$
was given by Molev, see \cite[Chapter 4]{M}.
In the present article
we will use another parametrization, introduced in \cite{KN}.
In the next subsection we will establish a correspondence
between the two parametrizations.
Let us call an $\Y(\so_n)\ts$-module \textit{integrable}
if the action of the Lie algebra $\so_n\subset\Y(\so_n)$
on it integrates to an action of the complex Lie group $\mathrm{SO}_n\ts$.
When working with the algebra $\Y(\so_n)\ts$,
we will be considering the integrable representations~only.

 
\bigskip\noindent
\textbf{1.5.}
First recall the parametrization from \cite{M}. 
Write $n=2l$ or $n=2l+1$ depending on whether $n$ is even or odd. If
$\g_n=\sp_n$ then $n$ has to be even.
Up to equivalence and similarity,
the irreducible finite-dimensional modules of the algebra $\Y(\sp_n)$ 
are parametrized by sequences of $l$
monic polynomials $Q_1(x)\lcd Q_l(x)$ with complex coefficients,
where the last polynomial $Q_l(x)$ is even. Further,
if $n$ is even then the integrable irreducible finite-dimensional 
$\Y(\so_n)\ts$-modules are parametrized by the same
sequences of polynomials, and by
an extra parameter $\de\in\{+1,-1\}$ in the case
when $Q_l\ts(0)=0\ts$.

If $n$ is odd then $\g_n=\so_n\ts$.
The integrable irreducible finite-dimensional 
modules of the algebra $\Y(\so_n)$ with odd $n$
are parametrized by sequences of $l$
monic polynomials $Q_1(x)\lcd Q_l(x)$ with complex coefficients,
but without any further conditions on the polynomial $Q_l(x)\ts$.
For the trivial module
of the algebra $\Y(\g_n)$ with arbitrary $\g_n$
we always have $Q_1(x)=\ldots=Q_l(x)=1\ts$, and
there is no extra parameter $\de$ then,
even if $\g_n=\so_n$ and $n=2l\ts$.

Now let $\f_m=\sp_{2m}$ if $\g_n=\sp_n$,
and let $\f_m=\so_{2m}$ if $\g_n=\so_n$.
Let $\h_m$ be a Cartan subalgebra of the reductive
Lie algebra $\f_m\ts$. The Weyl group of $\sp_{2m}$
is isomorphic to the hyperoctahedral group 
$H_m=S_m\ltimes\ZZ_2^m\ts$. 
The Weyl group of $\so_{2m}$
is isomorphic to a subgroup of $H_m$ of index two.
In the latter case the action of 
the Weyl group on $\h_m$ can be extended by
a diagram automorphism of $\mathrm{SO}_{2m}$ of order two,
so that the extended Weyl group is still isomorphic to $H_m\ts$.
Instead of 
$Q_1(x)\lcd Q_l(x)$ 
for any $\g_n$ we will use certain orbits in $\h_m^*\times\h_m^*$
under the diagonal shifted action of the group $H_m\ts$.

When working with the Lie algebra $\f_m$
it will be convenient to label the standard basis vectors
of $\CC^{2m}$ by the indices $-\ts m\lcd-1\com1\lcd m\ts$.
Let $a\com b$ be any pair of these indices. Let
$E_{ab}\in\gl_{2m}$ be the corresponding matrix unit.
Choose the antisymmetric bilinear form on $\CC^{2m}$
so that the subalgebra $\f_m\subset\gl_{2m}$ 
preserving this bilinear form is spanned by the elements
$$
F_{ab}=E_{ab}-\sgn\ts(a\ts b)\ts\cdot E_{-b,-a}
\quad\text{or}\quad
F_{ab}=E_{ab}-E_{-b,-a}
$$
for $\f_m=\sp_{2m}$ or $\f_m=\so_{2m}$ respectively.
Choose the Cartan subalgebra $\h_m$ of $\f_m$ 
with the basis $F_{11}\lcd F_{mm}\ts$.
For any weight $\la\in\h_m^*$ define the sequence
$\la_1\lcd\la_m$ of its {labels} by 
$\la_a=\la\ts(F_{aa})$ 
for $a=1\lcd m\ts$. For the half-sum $\rho\in\h_m^*$
of positive roots of $\f_m$ we have 
$\rho_a=-\ts a$ if $\f_m=\sp_{2m}\ts$,
and $\rho_a=1-a$ if $\f_m=\so_{2m}\ts$.
Here the positive roots are chosen as in \cite[Subsection 4.1]{KN}.

The group $H_m$ acts on the Cartan subalgebra $\h_m$ by permuting the 
basis vectors chosen above, and by multiplying any
of them by $-1\ts$. The group $H_m$ also acts on the dual space $\h_m^*$.
The shifted action of any element $w\in H_m$
is given by the universal formula \eqref{shiftact}.
Let $\ka\in\h_m^*$ be the weight of $\f_m$ such that 
$\ka_a=n/2$ for $a=1\lcd m\ts$.
Instead of the sequences of polynomials $Q_1(x)\lcd Q_l(x)$
for our parametrization of irreducible $\Y(\g_n)\ts$-modules we 
use the orbits of the pairs $(\la\com\mu)\in\h_m^*\times\h_m^*$
such that all labels of the weight $\nu=\la-\mu+\ka$ belong to the set
$\{\tts1\lcd n-1\}\ts$. Note that the definitions of 
$\nu$ here and in the case of the Yangian $\Y(\gl_n)$ are different.

Given any such a pair $(\la\com\mu)$ define 
a sequence $Q_1(x)\lcd Q_l(x)$ of monic polynomials as follows.
For any $\g_n$ and each $i=1\lcd l$ put
\begin{equation}
\label{qi}
Q_i(x)=
\prod_{\nu_a=i}\,(x+\mu_a+\rho_a)
\,\ts\cdot\hspace{-6pt}
\prod_{\nu_a=n-i}\!\!(x-\mu_a-\rho_a)\ts.
\end{equation}
Here $\mu_1\lcd\mu_m$ and $\nu_1\lcd\nu_m$ are
the labels of the weights $\mu$ and $\nu$ of $\f_m\ts$.
Note that if $n$ is even then $l=n-l\ts$, so that
$Q_l(x)$ is an even polynomial~then.

The simultaneous shifted action of the subgroup $S_m\subset H_m$ on
the weights $\la$ and $\mu$ gives a permutational action of $S_m$
on the labels of $\mu+\rho$ and of $\nu\ts$. 
Further, for any $a=1\lcd m$
multiplying the basis vector $F_{aa}\in\h_m$ by $-1$
results in changing the labels $\mu_a+\rho_a$ and $\nu_a$
to respectively $-\mu_a-\rho_a$ and $n-\nu_a\ts$.
Therefore each of polynomials \eqref{qi}
depends only on the $H_m\ts$-orbit of the pair 
$(\la\com\mu)\in\h_m^*\times\h_m^*\ts$.
Moreover, the orbit is determined by these polynomials uniquely.
Furthermore, any sequence of monic polynomials $Q_1(x)\lcd Q_l(x)$
of total degree $m$ with complex coefficients arises in this way,
provided that for an even $n$ the last
polynomial $Q_l(x)$ is also even.

In the next subsection, to each of these orbits
we will attach an irreducible $\Y(\g_n)\ts$-module,
unless $\g_n=\so_n$ with even $n$ and $\mu_a+\rho_a=0$ for
some $a\ts$. In the latter case, to such an orbit
we will attach two irreducible $\Y(\g_n)\ts$-modules,
not equivalent to each other.
For these two, we will have $\de=1$ and $\de=-1$.

In any case, the polynomials $Q_1(x)\lcd Q_l(x)$  
of the attached modules will be given by \eqref{qi}.
For the definition of $Q_1(x)\lcd Q_l(x)$
for any irreducible finite-dimensional $\Y(\g_n)\ts$-module
see \cite[Subsections 5.3, 5.4, 5.5]{KN}.

 
\bigskip\noindent
\textbf{1.6.}
For $k=0\com1\lcd n$ let us denote by $\Ph^{\ts-k}_{\tts t}$ 
the $\Y(\gl_n)\ts$-module  
obtained by pulling the $\Y(\gl_n)\ts$-module $\Ph^{\ts k}_t$ back through
the automorphism \eqref{transauto}. Note that due to Lemma \ref{1.15}
the $\Y(\gl_n)\ts$-module
$\Ph^{\ts-k}_{\tts t}$ is similar to the module
$\Ph^{\ts n-k}_{1-t}$.

Take any $(\la\com\mu)\in\h_m^*\times\h_m^*$
such that all labels of the weight $\nu=\la-\mu+\ka$ belong to
the set $\{\tts1\lcd n-1\}\ts$. Let the weights 
$\la$ and $\mu$ vary so that $\nu$ is fixed.
For the Cartan subalgebra $\h_m$ of $\f_m$
the weight $\mu$ is \textit{generic\/} 
if $\mu_a-\mu_b\notin\ZZ$ and $\mu_a+\mu_b\notin\ZZ$
whenever $a\neq b\ts$, and $2\ts\mu_a\notin\ZZ$ for any $a\ts$.

Consider the tensor product of $\Y(\gl_n)\ts$-modules \eqref{basprod1}
where $\mu\com\nu$ and $\rho$ are now weights of $\f_m\ts$, not of
$\gl_m$ as before. 
It is known that the restriction of the $\Y(\gl_n)\ts$-module 
\eqref{basprod1} to the subalgebra $\Y(\g_n)\subset\Y(\gl_n)$
is irreducible if (but not only if) 
the weight $\mu$ of $\f_m$ is generic, 
see \cite[Theorems 5.3, 5.4, 5.5]{KN}.
Moreover, if (but not only if) $\mu$ is generic then 
all the $\Y(\g_n)\ts$-modules obtained from \eqref{basprod1}
by permuting the tensor factors and by replacing any tensor factor
$\Ph^{\ts k}_t$ by $\Ph^{\ts-k}_t$,
are equivalent to each other. In particular, for each generic 
weight $\mu$ of $\f_m$
there is a unique non-zero $\Y(\g_n)\ts$-intertwining operator
\begin{equation}
\label{intoper2}
\Ph^{\,\nu_1}_{\mu_1+\rho_1+\frac12}
\ot\ts\ldots\ts\ot\ts
\Ph^{\,\nu_m}_{\mu_m+\rho_m+\frac12}
\to\,
\Ph^{\,-\nu_1}_{\mu_1+\rho_1+\frac12}
\ot\ts\ldots\ts\ot\ts
\Ph^{\,-\nu_m}_{\mu_m+\rho_m+\frac12}
\end{equation}
corresponding to the element of the group $H_m$  
of the maximal length. We will
denote this operator by $J(\mu)\ts$;
it is unique up to a multiplier from $\CC\,\backslash\ts\{0\}\ts$.

For all generic weights $\mu$ of $\f_m\ts$,
the source and the target vector spaces of the operators
$J(\mu)$ are the same tensor product
of exterior powers of $\CC^n$
\begin{equation*}
\label{lamprod}
\La^{\nu_1}(\CC^n)
\ot\ts\ldots\ts\ot\ts
\La^{\nu_m}(\CC^n)\ts.
\end{equation*}
Hence $J(\mu)$ is a function of
$\mu$ taking values in the space of linear operators on this tensor product.
The multipliers from $\CC\,\backslash\{0\}$
can be chosen so that $J(\mu)$ becomes a rational function of $\mu\ts$,
see Section 3 for a particular choice. Any such a choice
allows to determine intertwining operators \eqref{intoper2}
for those non-generic weights $\mu$ where the function $J(\mu)$ is regular.
The next theorem summarizes the results of 
\cite[Subsections 5.3, 5.4, 5.5]{KN}.
It will be generalized in Subsection 3.7 of the present article.

For the Lie algebra $\f_m$ 
the weight $\la$ is \textit{dominant\/} if 
$\la_a-\la_b\neq-\,1\com-\,2,\,\ldots$ and $\la_a+\la_b\neq\,1,2,\,\ldots$ 
for all $a<b\ts$, with an extra condition that $\la_a\neq\,1,2,\,\ldots$
in the case $\f_m=\sp_{2m}\ts$.
The pair $(\la\com\mu)\in\h_m^*\times\h_m^*$
is called \textit{good} 
if the weight $\la+\rho$ is dominant and
\begin{align}
\label{con1}
\nu_a\ge\nu_b
&\quad\text{whenever}\quad
\la_a-\la_b+\rho_a-\rho_b=0
\quad\text{and}\quad
a<b\ts,
\\
\label{con2}
\nu_a+\nu_b\le n\ 
&\quad\text{whenever}\quad
\la_a+\la_b+\rho_a+\rho_b=0
\quad\text{and}\quad
a<b\ts,
\end{align}
with an extra condition that in the case $\f_m=\sp_{2m}$
\begin{equation}
\label{con3}
2\ts\nu_a\le n \ \ 
\quad\text{whenever}\quad
\la_a+\rho_a=0\ts.
\end{equation}
The orbit of any $(\la\com\mu)\in\h_m^*\times\h_m^*$
under the shifted action of $H_m$ on $\h_m^*$
does contain a good pair. Here we assume only that
$\nu_1\lcd\nu_m\in\{\tts1\lcd n-1\}\ts$. 
Note that for $\f_m=\so_{2m}$ a good pair 
is contained already in the orbit of any $(\la\com\mu)$
under the shifted action of the Weyl group,
which is a subgroup of $H_m$ of index two. 
However, this fact will not be used in the present article.


\begin{theorem*}
\label{1.2}
For the fixed $\nu=\la-\mu+\ka$
the multipliers from $\CC\,\backslash\{0\}$ of the 
operators \eqref{intoper2} for all generic weights $\mu$ can be chosen so 
that the rational function $J(\mu)$ is regular and 
non-zero whenever the pair $(\la\com\mu)$ is good.
Then for any good pair $(\la\com\mu)$
the quotient of \eqref{basprod1}
by the kernel of the operator $J(\mu)$
is an irreducible $\Y(\g_n)\ts$-module,
unless\/ $\g_n=\so_n$ with even $n$ and $\mu_a+\rho_a=0$
for some index\/ $a\ts$. In the latter case 
the quotient of \eqref{basprod1}
by the kernel of\/ $J(\mu)$
is a direct sum of two irreducible 
$\Y(\g_n)\ts$-modules, not equivalent to each other. 
For any\/ $\g_n$ the polynomials $Q_1(x)\lcd Q_l(x)$  
of the irreducible $\Y(\g_n)\ts$-modules occurring as above
are given by \eqref{qi}. 
\end{theorem*}


\section*{\normalsize\bf 2. Intertwining operators}
\label{S2}
\setcounter{section}{2}
\setcounter{theorem*}{0}
\setcounter{equation}{0}


\textbf{2.1.}
In this subsection we develop the formalism of $R$-matrices\ts;
it will be used to produce explicit formulas for
intertwining operators \eqref{intoper1} and \eqref{intoper2}
over $\Y(\gl_n)$ and $\Y(\g_n)$ respectively.   
Let $P$ denote the linear operator on $(\CC^n)^{\tts\ot\tts2}$
exchanging the two tensor factors. The \textit{Yang $R$-matrix}
is the rational function of a variable $x$
$$
R(x)=1-P\ts x^{-1}
$$
taking values in 
$$
\End(\CC^n)^{\tts\ot\tts2}=(\ts\End\CC^n)^{\tts\ot\tts2}\ts.
$$
It satisfies the \textit{Yang-Baxter equation} in 
$(\ts\End\CC^n)^{\tts\ot\tts3}$
\begin{equation}
\label{ybe}
R_{12}(x)\,R_{13}(x+y)\,R_{23}(y)=R_{23}(y)\,R_{13}(x+y)\,R_{12}(x)\,.
\end{equation}
As usual, the subscripts in \eqref{ybe} indicate
different embeddings of the algebra $(\ts\End\CC^n)^{\tts\ot\tts2}$
to $(\ts\End\CC^n)^{\tts\ot\tts3}$ so that
$R_{12}(x)=R(x)\ot1$ and $R_{23}(y)=1\ot R(y)\ts$.

We have
$$
P=\sum_{i,j=1}^n\,E_{ij}\ot E_{ji}\ts.
$$
Denote
\begin{equation}
\label{pb}
\Pb=\sum_{i,j=1}^n\,E_{ij}\ot E_{ij}
\quad\ \text{and}\ \quad
\Pt=\sum_{i,j=1}^n\,\widetilde{E}_{ij}\ot E_{ji}\ts.
\end{equation}
Then put 
$$
\Rb(x)=1-\Pb\ts x^{-1}
\quad\ \text{and}\ \quad
\Rt(x)=1-\Pt\ts x^{-1}\ts.
$$
The values of the function $\,\Rb(x)$
are obtained from those of $R(x)$ by applying the 
matrix transposition to the first tensor factor of
$(\ts\End\CC^n)^{\tts\ot\tts2}\ts$.
The values of the function $\Rt(x)$
are obtained from those of $R(x)$ by 
applying to the first tensor factor of
$(\ts\End\CC^n)^{\tts\ot\tts2}$
the conjugation with respect to the form 
$\langle\ ,\,\rangle\ts$.
Now \eqref{ybe} and the relation
$P\ts R(x)\ts P=R(x)$ imply that
\begin{align}
\label{tybe}
\Rb_{12}(x+y)\,\ts\Rb_{13}(x)\,R_{23}(y)&=
R_{23}(y)\,\ts\Rb_{13}(x)\,\ts\Rb_{12}(x+y)\,,
\\
\label{cybe}
\Rt_{12}(x+y)\,\Rt_{13}(x)\,R_{23}(y)&=
R_{23}(y)\,\Rt_{13}(x)\,\Rt_{12}(x+y)\,.
\end{align}

Finally, denote
$$
\widehat{P}=\sum_{i,j=1}^n\,\widetilde{E}_{ij}\ot E_{ij}
$$
and put
$$
\Rh(x)=1-\widehat{P}\ts x^{-1}\ts.
$$
The values of the function $\Rh(x)$
are obtained from those of $\Rb(x)$ by 
applying to the first tensor factor of
$(\ts\End\CC^n)^{\tts\ot\tts2}$
the conjugation with respect to $\langle\ ,\,\rangle\ts$.
Therefore the relation \eqref{tybe} implies that
\begin{equation}
\label{hybe}
\Rh_{12}(x)\,\Rh_{13}(x+y)\,R_{23}(y)=
R_{23}(y)\,\Rh_{13}(x+y)\,\Rh_{12}(x)\,.
\end{equation} 
Note that $P\,\Rb(x)\ts P=\Rb(x)$ and $P\ts\Rt(x)\ts P=\Rt(x)\ts$.
Further, due to 
$\Pb^2=n\,\Pb$ we have
$\Rb(x)\,\ts\Rb(n-x)=1\ts$. By using the latter relation,
\eqref{tybe} also implies that
\begin{equation}
\label{bybe}
\textrm{$
\Rb_{12}(n-x-y)\,\Rh_{13}(x)\,\Rt_{23}(y)=
\Rt_{23}(y)\,\Rh_{13}(x)\,\ts\Rb_{12}(n-x-y)\,.
$}
\end{equation}

Observe that $R(1)=1-P\ts$. A direct calculation now shows that
\begin{equation}
\label{nopole}
R_{12}(x)\,R_{13}(x+1)\,R_{23}(1)=
\bigl(\ts1-(\ts P_{12}+P_{13})\,x^{-1}\bigr)(\ts1-P_{23})\,.
\end{equation}
In particular, when $y=1$, 
the rational function of $x$ at either side of 
\eqref{ybe} has no pole at $x=-1\ts$. Similarly, when $y=1\ts$,
the rational function of $x$ at either side of \eqref{tybe}
or \eqref{cybe} has no pole at $x=-1\ts$. Moreover,
\begin{equation}
\label{nopolet}
\Rt_{12}(x+1)\,\Rt_{13}(x)\,R_{23}(1)=
\bigl(\ts1-(\ts\Pt_{12}+\Pt_{13})\,x^{-1}\bigr)(\ts1-P_{23})\,.
\end{equation}


\smallskip\noindent
\textbf{2.2.}
Now consider the representation of $\Y(\gl_n)$
obtained by pulling the defining representation of $\gl_n$
back through the homomorphism \eqref{eval}, and then back through the
automorphism \eqref{tauz} of $\Y(\gl_n)\ts$.
The resulting $\Y(\gl_n)\ts$-module has been denoted by
$\Ph^{\tts 1}_{\tts t}\ts$. Note that under this representation
\begin{equation}
\label{teval}
T(x)\mapsto\Rb(\ts t-x)\ts.
\end{equation}
Here on the left we regard $n\times n$ matrices with entries
from the algebra $\Y(\gl_n)$ as elements of
$\End\CC^n\ot\Y(\gl_n)\ts$;
we will always do so in this and in the next section.
Our explicit formulas for intertwining operators
\eqref{intoper1} and \eqref{intoper2} are based 
on the following simple and well known lemma,
first appeared in \cite{KRS}. 

\begin{lemma*}
\label{2.1}
For any $k=1\lcd n$ and\/ $t\in\CC$ the $\Y(\gl_n)\ts$-module 
$\Ph^{\ts k}_t$ appears as the submodule of 
\begin{equation}
\label{kext}
\Ph^{\tts1}_{\tts t+k-1}
\ot\ts\ldots\ts\ot\ts
\Ph^{\tts1}_{\tts t+1}
\ot
\Ph^{\tts1}_{\tts t}
\end{equation}
with underlying subspace\/ 
\begin{equation}
\label{lak}
\La^k(\CC^n)\subset(\CC^n)^{\tts\ot\tts k}\ts.
\end{equation}
\end{lemma*}

\begin{proof}
First consider the standard action of the Lie algebra $\gl_n$
on the vector space $(\CC^n)^{\tts\ot\tts k}\ts$. 
Turn this vector space into an $\Y(\gl_n)\ts$-module,
by pulling the action of $\gl_n$
back through the homomorphism \eqref{eval} and then back through the
automorphism \eqref{tauz} of $\Y(\gl_n)\ts$.
Under the resulting representation of $\Y(\gl_n)$
\begin{equation}
\label{ppp}
T(x)
\,\mapsto\,
1+(\,\Pb_{01}+\ldots+\Pb_{0,k-1}+\Pb_{0k})\,(\ts x-t)^{-1}\ts,
\end{equation}
see \eqref{pb}.
Here we use the subscripts $0\lcd k-1\com k$ rather than
$1\lcd k\com k+1$ to label the tensor factors of
$(\ts\End\CC^n)^{\tts\ot\tts(k+1)}\ts$.
The $\Y(\gl_n)\ts$-module $\Ph^{\ts k}_t$
is defined by restricting the above described action 
of $\Y(\gl_n)$ to the subspace \eqref{lak}.

Next consider the action of $\Y(\gl_n)$ on \eqref{kext}.
By \eqref{1.33} and \eqref{teval}, then 
$$
T(x)
\,\mapsto\,
\Rb_{\ts01}(\ts t+k-1-x)
\,\ldots\,
\Rb_{\ts0,k-1}(\ts t+1-x)\,\ts\Rb_{\ts0k}(\ts t-x)\ts. 
$$
The latter product is a rational function of $x\ts$, 
valued in $(\ts\End\CC^n)^{\tts\ot\tts(k+1)}$.
It tends to $1$ when $x\to\infty\ts$, 
and has a simple pole at $x=t$ with the residue
$$
P_{01}+\ldots+P_{0,k-1}+P_{0k}\ts.
$$
Further, by an observation made after \eqref{nopole}
for any $i=1\lcd k-1$ the product
$$
\Rb_{\ts01}(\ts t+k-1-x)
\,\ldots\,
\Rb_{\ts0,k-1}(\ts t+1-x)\,\ts\Rb_{\ts0k}(\ts t-x)
\ts(\ts1-P_{i,i+1})
$$
has no pole at $x=t+k-i\ts$.
Hence the restriction of
the action of $\Y(\gl_n)$ on \eqref{kext}
to the subspace \eqref{lak}
coincides with the restriction~of~the action of $\Y(\gl_n)$
on the vector space
$(\CC^n)^{\tts\ot\tts k}$ described by the assignment \eqref{ppp}.
\qed
\end{proof}

The assignment \eqref{transauto} defines
a coalgebra anti-automorphism of $\Y(\gl_n)\ts$.
This immediately implies another lemma,
to be used when working with $\Y(\g_n)\ts$.

\begin{lemma*}
\label{2.2}
The $\Y(\gl_n)\ts$-module obtained by pulling
\eqref{kext} back through the automorphism \eqref{transauto}
is equivalent to 
\begin{equation*}
\Ph^{\ts-1}_{\tts t}
\ot
\Ph^{\ts-1}_{\tts t+1}
\ot\ts\ldots\ts\ot\ts
\Ph^{\ts-1}_{\tts t+k-1}\ts.
\end{equation*}
The linear operator on $(\CC^n)^{\tts\ot\tts k}$
reversing the order of tensor factors
intertwines the two equivalent modules.
\end{lemma*}

For $k=0\com1\lcd n$ denote respectively by
$\Phd^{\ts k}_t$ and $\Phd^{\ts-k}_{\tts t}$
the $\Y(\gl_n)\ts$-modules 
obtained by pulling $\Ph^{\ts k}_t$ and $\Ph^{\ts-k}_{\tts t}$
back through the automorphism \eqref{fut} where 
\begin{equation*}
\label{fxt}
f(x)=\frac{x-t}{x-t+1}\,\ts.
\end{equation*}

\begin{lemma*}
\label{1.15}
The\/ $\Y(\gl_n)\ts$-modules $\Ph^{\ts-k}_{\tts t}$ and\/
$\Phd^{\ts n-k}_{1-t}$ are mutually equivalent.
\end{lemma*}

\begin{proof}
Let $e_1\lcd e_n$ be the standard basis vectors of $\CC^n\ts$.
Denote by $x_i$ the operator of left multiplication
by the element $e_i$ in the exterior algebra $\La\ts(\CC^n)\ts$.
Let $\d_i$ be the corresponding operator of left derivation
on $\La(\CC^n)\ts$.
Note that 
\begin{equation*}
\label{ed}
x_i\ts\d_j+\d_j\ts x_i=\de_{ij}\,.
\end{equation*}
For each $k=0\com1\lcd n$ the action of the algebra
$\Y(\gl_n)$ on its module $\Ph^{\ts k}_t$
is defined by the homomorphism $\Y(\gl_n)\to\End(\La\ts(\CC^n))$ which maps
\begin{equation}
\label{phikat}
T_{ij}(x)\,\mapsto\,\de_{ij}+\frac{\,x_i\,\d_j\ts}{x-t}\ts\,.
\end{equation}
 
It suffices to prove Lemma \ref{1.15} for any choices
of the symmetric and of the alternating form    
$\langle\ ,\,\rangle$ on $\CC^n\ts$. For the proof only, 
choose the form as follows.
Put $\th_i=-1$ if $\g_n=\sp_n$ and $i>n/2\ts$; otherwise put $\th_i=1\ts$.
Set $\langle\ts e_i\com e_j\ts\rangle=\th_i\,\de_{\ts\bi\ts j}$
where we write $\bi=n-i+1$ for short. Then
\begin{equation}
\label{teij}
\widetilde{E}_{\ts ij}=\th_i\ts\th_j\,E_{\ts\bj\ts\bi}\,.
\end{equation}
Using \eqref{transauto} and \eqref{phikat},
the action of 
$\Y(\gl_n)$ on $\Ph^{\ts-k}_{\tts t}$
is then defined by 
\begin{equation}
\label{modone}
T_{ij}(x)\,\mapsto\,\de_{ij}-
\frac{\,\th_i\,\th_j\,x_{\ts\bj}\,\d_{\ts\bi}\ts}{x+t}\ts\,.
\end{equation}

On the other hand, the action of the algebra
$\Y(\gl_n)$ on its module $\ts\Phd^{\ts n-k}_{1-t}\ts$
is defined by the homomorphism 
$\Y(\gl_n)\to\End(\La\ts(\CC^n))$ which maps
\begin{equation}
\label{modtwo}
T_{ij}(x)\,\mapsto\,
\frac{x+t-1}{x+t}\,
\Bigl(\,\de_{ij}+\frac{\,x_i\,\d_j\ts}{x+t-1}\,\Bigr)
\,=\,
\de_{ij}-\frac{\d_j\ts x_i}{x+t}
\ts\,.
\end{equation}
By comparing the right hand sides of the assignments
\eqref{modone} and \eqref{modtwo},
the equivalence of the $\Y(\gl_n)\ts$-modules 
$\Phd^{\ts n-k}_{1-t}$
and
$\Ph^{\ts-k}_{\tts t}$ 
can be realized by the linear map of the underlying vector spaces
$\La^{n-k}(\CC^n)\to\La^k(\CC^n)\ts$:
$$
e_{i_1}
\wedge\ts\ldots\ts\wedge\ts
e_{i_{n-k}}
\mapsto
(\,\th_{j_1}\ts e_{\ts\bj_1}\ts)
\ts\wedge\ts\ldots\ts\wedge\ts
(\,\th_{j_k}\ts e_{\ts\bj_k}\ts)
$$
where 
$$
e_{i_1}
\wedge\ts\ldots\ts\wedge\ts
e_{i_{n-k}}
\wedge
e_{j_1}
\wedge\ts\ldots\ts\wedge\ts
e_{j_k}
\,=\,
e_1
\wedge\ts\ldots\ts\wedge\ts
e_n\,\ts.
$$
Indeed, for any indices $i\com j=1\lcd n$ this map 
intertwines the operators $\d_j\ts x_i$ and 
$\th_i\,\th_j\,x_{\ts\bj}\,\d_{\ts\bi}$ on the vector spaces
$\La^{n-k}(\CC^n)$ and $\La^k(\CC^n)$ respectively.
\qed
\end{proof}


\noindent
\textbf{2.3.}
Once again take any weights $\la$ and $\mu$ of $\gl_m$
such that all the labels of the weight $\nu=\la-\mu$ belong to
the set $\{\tts1\lcd n-1\}\ts$. Put $N=\nu_1+\ldots+\nu_m$  
and split the sequence $1\lcd N$ to the
consecutive segments of lengths $\nu_1\lcd\nu_m\ts$.
Hence the $a\ts$th segment is the sequence of numbers 
\begin{equation}
\label{pnu}
p=\nu_1+\ldots+\nu_{a-1}+i 
\quad\text{where}\quad
i=1\lcd\nu_a\ts.
\end{equation}
Then put
\begin{equation}
\label{xp}
\textstyle
x_p=\mu_a+\rho_a+\frac12+\nu_a-i\ts.
\end{equation}
Let
$$
P_\nu:\ts
(\CC^n)^{\tts\ot\ts\nu_1}
\ot\ts\ldots\ts\ot
(\CC^n)^{\tts\ot\ts\nu_m}
\ts\to\ts
(\CC^n)^{\tts\ot\ts\nu_m}
\ot\ts\ldots\ts\ot
(\CC^n)^{\tts\ot\ts\nu_1}
$$
be the linear operator on $(\CC^n)^{\tts\ot\tts N}$
reversing the order of the tensor factors $\CC^n$
by segments of lengths
$\nu_1\lcd\nu_m$ in their sequence. Within any segment, the order 
of tensor factors is not changed. Then
$P_\nu=P$ if $m=2$ and $\nu_1=\nu_2=1\ts$.

Now let the weight $\mu\in\t_m^*$ vary while 
the weight $\nu$ is fixed.  
Let $1'\lcd N'$ be the sequence obtained from $1\lcd N$
by reversing the order of the terms by the segments of lengths
$\nu_1\lcd\nu_m$ introduced above. Within every segment,
the order of terms is not changed.
Let $1''\lcd N\tts''$ be the sequence obtained from $1\lcd N$
by reversing the order of terms within the segments.
The order of the segments themselves is not changed now. 
Take the ordered product
\begin{equation}
\label{bmu}
B(\mu)=\prod_{(p,q)}^{\longrightarrow}
\,R_{\ts pq}\ts(\tts x_q-x_p\tts)
\end{equation}
where $p<q$ and they belong to different segments of the sequence
$1\lcd N\ts$. Here the pair 
$(p\com q)$ precedes $(r\ts,s)$ 
if $p<r$ or if $p=r$ and $q$ precedes $s$ in the sequence $1''\lcd N''$. 
Note that $B(\mu)$ is a rational function of $\mu$
without poles at generic weights of $\gl_m\ts$.

\begin{proposition*}
\label{2.3}
Suppose that the weight $\mu$ of\/ $\gl_m$ is generic.
Then $P_\nu\ts B(\mu)$ is an intertwining operator of\/ 
$\Y(\gl_n)\ts$-modules
\begin{equation}
\label{pa}
\Ph^{\tts1}_{\tts x_1}
\ot\ldots\ts\ot\ts
\Ph^{\tts1}_{\tts x_N}
\to\,
\Ph^{\tts1}_{\tts x_{1'}}
\ot\ldots\ts\ot\ts
\Ph^{\tts1}_{\tts x_{N'}}
\ts.
\end{equation}
\end{proposition*}

\begin{proof}
Under the action of $\Y(\gl_n)$ on the source module in \eqref{pa},
\begin{equation}
\label{x}
T(x)
\,\mapsto\,
\Rb_{\ts01}(x_1-x)
\,\ldots\,\ts
\Rb_{\ts0N}(\ts x_N-x)\ts.
\end{equation}
Like in the proof of Lemma \ref{2.1},
here we use the subscripts $0\com1\lcd N$ rather than
$1\com2\lcd N+1$ to label the tensor factors of
$(\ts\End\CC^n)^{\tts\ot\tts(N+1)}\ts$.
Similarly,
under the action of $\Y(\gl_n)$ on the target module in \eqref{pa},
\begin{equation}
\label{x'}
T(x)
\,\mapsto\,
\Rb_{\ts01}(x_{1'}-x)
\,\ldots\,\ts
\Rb_{\ts0N}(\ts x_{N'}-x)\ts.
\end{equation}
Denote by $X$ and $X'$ the right hand sides of 
the assignments \eqref{x} and \eqref{x'} respectively.
By using the relation \eqref{tybe} repeatedly, we get
$$
P_\nu\ts B(\mu)\,X=
P_\nu\,
\Rb_{\ts01'}(x_{1'}-x)
\,\ldots\,\ts
\Rb_{\ts0N'}(\ts x_{N'}-x)
\,B(\mu)=
X'\ts P_\nu\ts B(\mu)\ts.
$$
The equality of the left and of the right hand sides here
proves the claim.
\qed
\end{proof}


\noindent
\textbf{2.4.}
Let $A_k$ be the operator of antisymmetrization
on $(\CC^n)^{\tts\ot\tts k}\ts$, normalized so that $A_k^2=A_k\ts$.  
The subspace \eqref{lak} is the image of $A_k\ts$.
The ordered product
\begin{equation}
\label{ak}
\prod_{(i,j)}^{\longrightarrow}
\,R_{\ts ij}\ts(j-i)
\,=\,k\tts!\,A_k
\end{equation}
where $1\le i<j\le k$ and the pairs $(i\com j)$ 
are ordered lexicographically.
The formula \eqref{ak} has appeared in \cite{KRS}
but was already known to Jucys \cite{J}. 

Let $e_1\lcd e_n$ be the standard basis vectors of $\CC^n\ts$. 
For each $k=1\lcd n$ consider the vector
$\ph_k=e_1\wedge\ldots\wedge e_k\in\La^k(\CC^n)\ts$.
Using the embedding \eqref{lak},
\begin{equation*}
\label{phik}
\ph_k=A_k\ts(\tts e_1\ot\ts\ldots\ts\ot e_k\tts)\ts.
\end{equation*}
The next proposition is known\ts;
see \cite[Theorem 2]{N} for a more general result.
It is still instructive to give a proof here,
as it will be used later on.

\begin{proposition*}
\label{2.4}
For any generic weight 
$\mu$ of\/ $\gl_m$ 
the vector 
$\ph_{\nu_1}\ot\ts\ldots\ts\ot\ts\ph_{\nu_m}$
is an eigenvector of the operator $B(\mu)$ on
$(\CC^n)^{\tts\ot\tts N}$ with the eigenvalue
\begin{equation}
\label{beigen}
\prod_{1\le a<b\le m}\ 
\left\{
\begin{array}{ll}
\displaystyle
\ \frac{\,\la_a-\la_b+\rho_a-\rho_b+\nu_b}
{\,\mu_a-\mu_b+\rho_a-\rho_b}
&\quad\textrm{if}\quad\ \nu_a\le\nu_b\,;
\\[16pt]
\displaystyle
\ \frac{\,\la_a-\la_b+\rho_a-\rho_b+\nu_b}{\,\la_a-\la_b+\rho_a-\rho_b}
&\quad\textrm{if}\quad\ \nu_a\ge\nu_b\,.
\end{array}
\right.   
\end{equation}
\end{proposition*}

\begin{proof}
This proposition immediately follows from its particular case
of $m=2\ts$. Let us consider this case only. Then we have
\begin{gather}
\nonumber
B(\mu)\,(\ts\ph_{\nu_1}\ot\ph_{\nu_2}\tts)
\,=\!
\prod_{i=1\lcd\nu_1}^{\longrightarrow}
\!\Bigl(\,
\prod_{j=1,\ldots,\nu_2}^{\longleftarrow}
R_{\ts i,\nu_1+j}\ts(\tts x_{\nu_1+j}-x_i\tts)
\Bigr)\,\,\times
\\[6pt]
\nonumber
(\ts A_{\nu_1}\ot A_{\nu_2}\tts)\,
(\tts e_1\ot\ts\ldots\ts\ot e_{\nu_1}
\ot 
e_1\ot\ts\ldots\ts\ot e_{\nu_2}\tts)\,=\,
\\[6pt]
\label{aar}
(\ts A_{\nu_1}\ot A_{\nu_2}\tts)
\prod_{i=1,\ldots,\nu_1}^{\longleftarrow}
\Bigl(\,
\prod_{j=1,\ldots,\nu_2}^{\longrightarrow}
R_{\ts i,\nu_1+j}\ts(\tts x_{\nu_1+j}-x_i\tts)
\Bigr)
\\[6pt]
\label{eeee}
(\tts e_1\ot\ts\ldots\ts\ot e_{\nu_1}
\ot 
e_1\ot\ts\ldots\ts\ot e_{\nu_2}\tts)
\end{gather}
where the last equality is obtained by using the formula \eqref{ak}
and by applying 
\eqref{ybe} repeatedly. 
The reversed arrow over the product symbol indicates
that the factors corresponding to the running index 
are arranged from right to left.

First suppose that $\nu_1\le\nu_2\ts$. 
Arguing like in the proof of Lemma \ref{2.1} 
we~can always rewrite the product displayed in line \eqref{aar} as
$$
(\ts A_{\nu_1}\ot A_{\nu_2}\tts)
\prod_{i=1,\ldots,\nu_1}^{\longleftarrow}
\Bigl(\ts1-\ts
(x_{\nu_1+\nu_2}-x_i)^{-1}\,
\sum_{j=1}^{\nu_2}\,
P_{\ts i,\nu_1+j}
\,\Bigr)
$$ 
where any sum over $j=1\lcd\nu_2$ clearly
commutes with the operator $1\ot A_{\nu_2}$ on
$(\CC^n)^{\tts\ot\ts(\nu_1+\nu_2)}$. But for $\nu_1\le\nu_2$
the operators
$(1\ot A_{\nu_2}\tts)\,P_{\ts i,\nu_1+j}$
with $i\neq j$ annihilate the vector \eqref{eeee}
while the operators $P_{\ts i,\nu_1+i}$
do not change it. 
Hence applying the operator \eqref{aar} to \eqref{eeee}
gives the vector $\ph_{\nu_1}\ot\ph_{\nu_2}$ multiplied~by
$$
\prod_{i=1,\ldots,\nu_1}
\bigl(\ts1-\ts
(x_{\nu_1+\nu_2}-x_i)^{-1}
\ts\bigr)\,=\,
\frac{\,\la_1-\la_2+\rho_1-\rho_2+\nu_2}{\,\mu_1-\mu_2+\rho_1-\rho_2}\ .
$$

Next suppose that $\nu_1\ge\nu_2\ts$. 
We can always rewrite the product \eqref{aar} as
$$
(\ts A_{\nu_1}\ot A_{\nu_2}\tts)
\prod_{j=1,\ldots,\nu_2}^{\longrightarrow}
\Bigl(\,
\prod_{i=1,\ldots,\nu_1}^{\longleftarrow}
R_{\ts i,\nu_1+j}\ts(\tts x_{\nu_1+j}-x_i\tts)
\Bigr)\,.
$$
Arguing like in the proof of Lemma \ref{2.1} 
we~can rewrite the latter product as
$$
(\ts A_{\nu_1}\ot A_{\nu_2}\tts)
\prod_{j=1,\ldots,\nu_2}^{\longrightarrow}
\Bigl(\ts1-\ts
(x_{\nu_1+j}-x_1)^{-1}\,
\sum_{i=1}^{\nu_1}\,
P_{\ts i,\nu_1+j}
\,\Bigr)
$$
where any sum over $i=1\lcd\nu_1$ clearly
commutes with the operator $A_{\nu_1}\ot1$ on
$(\CC^n)^{\tts\ot\ts(\nu_1+\nu_2)}$. But for $\nu_1\ge\nu_2$
the operators
$(A_{\nu_1}\ot1\tts)\,P_{\ts i,\nu_1+j}$
with $i\neq j$ annihilate the vector \eqref{eeee}
while the operators $P_{\ts i,\nu_1+i}$
do not change it. 
Hence applying the operator \eqref{aar} to \eqref{eeee}
gives the vector $\ph_{\nu_1}\ot\ph_{\nu_2}$ multiplied~by
$$
\prod_{j=1,\ldots,\nu_2}
\bigl(\ts1-\ts
(x_{\nu_1+j}-x_1)^{-1}
\ts\bigr)\,=\,
\frac{\,\la_1-\la_2+\rho_1-\rho_2+\nu_2}{\,\la_1-\la_2+\rho_1-\rho_2}\ .
$$
This observation completes the proof of Proposition \ref{2.4}.
\qed
\end{proof}


\noindent
\textbf{2.5.}
Using the relations \eqref{ybe} and $R(1)=1-P$ repeatedly, 
one demonstrates that for any generic weight $\mu$ of $\gl_m$
the operator $B(\mu)$ preserves the subspace 
\begin{equation}
\label{subspace}
\La^{\nu_1}(\CC^n)
\ot\ts\ldots\ts\ot
\La^{\nu_m}(\CC^n)
\subset
(\CC^n)^{\tts\ot\tts N}\ts,
\end{equation}
see the proof of Proposition \ref{2.4} above.
Moreover, we have another proposition.

\begin{proposition*}
\label{2.45}
For any generic weight\/ $\mu$ of\/ $\gl_m$ 
the restriction of\/ 
$B(\mu)$ to the subspace \eqref{subspace} coincides with that of the operator
\begin{equation}
\label{baddit}
\prod_{1\le a<b\le m}^{\longrightarrow}\ 
\!\Bigl(\ 1\ +\   
\sum_{d>0}\ 
\sum_{\substack{i_1,\ldots,i_d\\j_1,\ldots,j_d}}\ \, 
\prod_{k=1}^d\ \ \ 
\frac{P_{\ts i_kj_k}}{\la_a-\la_b+\rho_a-\rho_b+\nu_b-k}
\ \Bigr)
\end{equation}
where we order the pairs $(a\com b)$ lexicographically
while $\,i_1\lcd i_d\,$ and $\,j_1\lcd j_d$ are pairwise 
distinct numbers respectively from the $a\,$th and $b\,$th 
segments of the sequence $1\lcd N$ taken so that
different are all corresponding sets of~$d$~pairs 
\begin{equation}
\label{dp}
(i_1\com j_1)\lcd(i_d\com j_d)\ts.
\end{equation}
\end{proposition*}

\begin{proof}
Note that as the indices $i_1\lcd i_d$ and $j_1\lcd j_d$ 
are pairwise distinct, we have $d\le\nu_a$ and
$d\le\nu_b$ for any non-zero summand in the brackets
in \eqref{baddit}.

The proposition immediately follows from its particular case
of $m=2\ts$. Let us consider this case only. Then \eqref{baddit}
is an operator on $(\CC^n)^{\tts\ot\ts(\nu_1+\nu_2)}$~equal~to
\begin{equation}
\label{atwo}
\sum_{d\ge0}\ 
\sum_{\substack{i_1,\ldots,i_d\\j_1,\ldots,j_d}}\ \, 
\prod_{k=1}^d\ \ \ 
\frac{P_{\ts i_kj_k}}{x_k-x_{\nu_1+\nu_2}}
\end{equation}
where $i_1\lcd i_d$ and $j_1\lcd j_d$ are pairwise 
distinct numbers taken respectively from 
$1\lcd \nu_1$ and $\nu_1+1\lcd\nu_1+\nu_2$ so that 
different are all the corresponding sets of $d$ pairs \eqref{dp}. 
Here we use the equalities
$$
\textstyle
x_{\nu_1+\nu_2}=\mu_2+\rho_2+\frac12
\ \quad\text{and}\ \quad
x_k=\mu_1+\rho_1+\frac12+\nu_1-k
$$
for any $k\le\nu_1$.
We also assume that $1$ is the only term in \eqref{atwo}
with $d=0\ts$. 
On the other hand, for $m=2$ by definition we have
\begin{equation}
\label{btwo}
B(\mu)
\,=\!
\prod_{i=1\lcd\nu_1}^{\longrightarrow}
\!\Bigl(\,
\prod_{j=1,\ldots,\nu_2}^{\longleftarrow}
R_{\ts i,\nu_1+j}\ts(\tts x_{\nu_1+j}-x_i\tts)
\Bigr)\,.
\end{equation} 
Let us relate two operators on the vector space 
$(\CC^n)^{\tts\ot\ts(\nu_1+\nu_2)}$
by the symbol~$\,\equiv$ if their actions coincide on 
the subspace
\begin{equation}
\label{subtwo}
\La^{\nu_1}(\CC^n)
\ot
\La^{\nu_2}(\CC^n)
\subset
(\CC^n)^{\tts\ot\ts(\nu_1+\nu_2)}\ts.
\end{equation}
We will establish the relation
$\,\equiv\,$ between \eqref{atwo} and \eqref{btwo}
by induction on~$\nu_1\ts$. 

If $\nu_1>1$ then we assume the latter relation holds for $\nu_1-1$ 
instead of $\nu_1\ts$; if $\nu_1=1$ then we are not making any assumption.
Arguing like in the proof of Lemma \ref{2.1}
and then using the induction assumption, \eqref{btwo}
is related by~$\,\equiv\,$~to
\begin{align*}
\prod_{i=1\lcd\nu_1-1}^{\longrightarrow}
\Bigl(\,\,
\prod_{j=1,\ldots,\nu_2}^{\longleftarrow}
R_{\ts i,\nu_1+j}\ts(\tts x_{\nu_1+j}-x_i\tts)
\Bigr)
&\times
\Bigl(\ 1\,+\,
\sum_{j=1}^{\nu_2}\,
\frac{P_{\ts\nu_1,\,\nu_1+j}}
{x_{\nu_1}-x_{\nu_1+\nu_2}}
\ \Bigr)\,
\\[6pt]
\equiv\ 
\Bigl(\ \ 
\sum_{d\ge0}\ 
\sum_{\substack{i_1,\ldots,i_d\\j_1,\ldots,j_d}}\ \, 
\prod_{k=1}^d\ \ \ 
\frac{P_{\ts i_kj_k}}{x_k-x_{\nu_1+\nu_2}}
\ \Bigr)
&\times
\Bigl(\ 1\,+\,
\sum_{j=1}^{\nu_2}\,
\frac{P_{\ts\nu_1,\,\nu_1+j}}
{x_{\nu_1}-x_{\nu_1+\nu_2}}
\ \Bigr)\,
\end{align*}
where $i_1\lcd i_d$ and $j_1\lcd j_d$ are distinct 
indices taken respectively from 
$1\lcd\nu_1-1$ and $\nu_1+1\lcd\nu_1+\nu_2\ts$. We assume that 
all corresponding sets of $d$ pairs \eqref{dp} are different. 
The right hand side of the last relation equals
\begin{gather} 
\label{lone}
\sum_{d\ge0}\ 
\sum_{\substack{i_1,\ldots,i_d\\j_1,\ldots,j_d}}\ \, 
\prod_{k=1}^d\ \ \ 
\frac{P_{\ts i_kj_k}}{x_k-x_{\nu_1+\nu_2}}
\ \ \ \,+
\\
\label{ltwo}
\sum_{d>0}\ \, 
\sum_{\substack{i_1,\ldots,i_d\\j_1,\ldots,j_d}}\ \, 
\sum_{l=1}^d\ \,
\Bigl(\ \
\prod_{k=1}^d\ \ 
\frac{P_{\ts i_kj_k}}{x_k-x_{\nu_1+\nu_2}}
\ \Bigr)\ 
\frac{P_{\ts\nu_1j_l}}
{x_{\nu_1}-x_{\nu_1+\nu_2}}
\ \ \ \,+
\\
\label{lthree}
\sum_{d\ge0}\ 
\sum_{\substack{i_1,\ldots,i_d\\j_1,\ldots,j_d}}\ 
\sum_j\, 
\Bigl(\ \
\prod_{k=1}^d\ \ 
\frac{P_{\ts i_kj_k}}{x_k-x_{\nu_1+\nu_2}}
\ \Bigr)\ 
\frac{P_{\ts\nu_1,\,\nu_1+j}}
{x_{\nu_1}-x_{\nu_1+\nu_2}}
\end{gather}
where the index 
$j$ is taken from $1\lcd\nu_2$ but $\,\nu_1+j\neq j_1\lcd j_d\,$ however.

Consider the sum displayed in the line \eqref{ltwo}. Here we have
\begin{gather*}
\Bigl(\
\prod_{k=1}^d\, 
P_{\ts i_kj_k}
\Bigr)\,
P_{\ts\nu_1j_l}
\,=\,
\Bigl(\
\prod_{k\neq l}\, 
P_{\ts i_kj_k}
\Bigr)\,
P_{\ts i_lj_l}\,
P_{\ts\nu_1j_l}
\,=
\\
\,\Bigl(\
\prod_{k\neq l}\, 
P_{\ts i_kj_k}
\Bigr)\,
P_{\ts\nu_1j_l}\,
P_{\ts i_k\nu_1}
\,\equiv\,
-\,\ts\Bigl(\
\prod_{k\neq l}\, 
P_{\ts i_kj_k}
\Bigr)\,
P_{\ts\nu_1j_l}
\end{gather*}
where 
the right hand side does not
involve the index $i_l\ts$. 
Now let us fix a number $j\in\{1\lcd\nu_2\}$
and take any set of $d$ pairs \eqref{dp}
such that one of the pairs contains the number $\nu_1+j\ts$.
Then this number has the form of $j_l$ for some index $l\ts$.
If the set of the other $d-1$ pairs $(i_k\com j_k)$ with $k\neq l$
is also fixed, then $i_l$ ranges over a set
of cardinality $\nu_1-d\ts$, namely over the fixed set 
$$
\{1\lcd\nu_1-1\}\setminus\{i_1\lcd i_{l-1}\com i_{l+1}\lcd i_d\}\ts.
$$

Now let us perform the summation over the indices $i_l\com j_l$ and $l$ 
in \eqref{ltwo} first of all the running indices. After that,
rename the running indices $i_{l+1}\lcd i_d$ and $j_{l+1}\lcd j_d$ 
respectively by $i_l\lcd i_{d-1}$ and $j_l\lcd j_{d-1}\ts$.
By the arguments given in the previous paragraph, the sum \eqref{ltwo}
gets related by $\,\equiv\,$ to the sum
\begin{equation}
\label{lfour}
\sum_{d>0}\, 
\sum_{\substack{i_1,\ldots,i_{d-1}\\j_1,\ldots,j_{d-1}}}\, 
\sum_j\ \ \,
\Bigl(\ \ \ 
\prod_{k=1}^{d-1}\ \ 
\frac{P_{\ts i_kj_k}}{x_k-x_{\nu_1+\nu_2}}
\ \Bigr)\ 
\frac{d-\nu_1}{x_d-x_{\nu_1+\nu_2}}
\ \,
\frac{P_{\ts\nu_1,\,\nu_1+j}}{x_{\nu_1}-x_{\nu_1+\nu_2}}\ \,
\end{equation}
where $i_1\lcd i_{d-1}$ and $j_1\lcd j_d$ are distinct 
indices taken respectively from $1\lcd\nu_1-1$ and
$\nu_1+1\lcd\nu_1+\nu_2$ so that different are  
all sets of $d-1$ pairs 
$$
(i_1\com j_1)\lcd(i_{d-1}\com j_{d-1})
$$
while the index $j$ is taken from $1\lcd\nu_2$ but 
$\,\nu_1+j\neq j_1\lcd j_{d-1}\,$ however.

Replace the running index $d\ge0$ in \eqref{lthree}
by $d-1$ where $d>0\ts$.~We~get
\begin{equation}
\label{lfive}
\sum_{d>0}\, 
\sum_{\substack{i_1,\ldots,i_{d-1}\\j_1,\ldots,j_{d-1}}}\, 
\sum_j\ \ \,
\Bigl(\ \ \ 
\prod_{k=1}^{d-1}\ \ 
\frac{P_{\ts i_kj_k}}{x_k-x_{\nu_1+\nu_2}}
\ \Bigr)\ 
\frac{P_{\ts\nu_1,\,\nu_1+j}}
{x_{\nu_1}-x_{\nu_1+\nu_2}}
\end{equation}
with the same assumptions on the running indices as in \eqref{lfour}.
By adding up together the sums \eqref{lfour} and \eqref{lfive}, we get
\begin{equation}
\label{lsix}
\sum_{d>0}\, 
\sum_{\substack{i_1,\ldots,i_{d-1}\\j_1,\ldots,j_{d-1}}}\, 
\sum_j\ \ \,
\Bigl(\ \ \ 
\prod_{k=1}^{d-1}\ \ 
\frac{P_{\ts i_kj_k}}{x_k-x_{\nu_1+\nu_2}}
\ \Bigr)\ 
\frac{P_{\ts\nu_1,\,\nu_1+j}}
{x_d-x_{\nu_1+\nu_2}}
\end{equation}
by the equality $\,x_d+d=x_{\nu_1}+\nu_1\ts$.
The sum of \eqref{lone} and \eqref{lsix} equals \eqref{atwo}.
Thus we have made the induction step.
\qed
\end{proof}

By Lemma \ref{2.1} and by Proposition \ref{2.3} the restriction 
of operator $P_\nu\ts B(\mu)$ to the subspace \eqref{subspace} 
is an $\Y(\gl_n)\ts$-intertwining operator \eqref{intoper1}.
Let $I(\mu)$ be this restriction 
divided by the rational function \eqref{beigen}. 
Then by Proposition~\ref{2.4}
\begin{equation}
\label{imu}
I(\mu):\ts
\ph_{\nu_1}\ot\ts\ldots\ts\ot\ph_{\nu_m}\mapsto\ts
\ph_{\nu_m}\ot\ts\ldots\ts\ot\ph_{\nu_1}\ts.
\end{equation}

\begin{theorem*}
\label{2.5}
For any fixed weight $\nu=\la-\mu$
the rational function $I(\mu)$ is regular 
at any point $\mu\in\t_m^*$ where the weight $\la+\rho$ is dominant.
\end{theorem*}

The operator valued rational function $I(\mu)$
does not vanish at any point $\mu\in\t_m^*$ due to the
normalization \eqref{imu}. 
The regularity of $I(\mu)$ was proved in \cite{KN}
for all $\mu$ where the pair $(\la\com\mu)$ is good\tts;
our Theorem \ref{2.5} is more general. 

In next two subsections, we give two proofs of Theorem \ref{2.5}.
Each of them provides an explicit formula for the operator
$I(\mu)$ whenever $\la+\rho$ is dominant. 
We give two proofs,
because the resulting formulas for $I(\mu)$ are quite different.

However, in both proofs we will use the following observation. 
Suppose that the weight $\la+\rho$ of $\gl_m$ is dominant. It means that 
\begin{equation}
\label{domone}
\la_a-\la_b+\rho_a-\rho_b\neq-1\com-2\com\,\ldots
\quad\text{whenever}\quad
a<b\,.
\end{equation}
Then
$$
\la_a-\la_b+\rho_a-\rho_b+\nu_b\neq0
$$
since $\nu_b$ is a positive integer. So 
the rational function \eqref{beigen} does not vanish then.
Moreover, then any factor of the product \eqref{beigen}
with $\nu_a\ge\nu_b$ has a simple pole as a function
of $\la_a-\la_b+\rho_a-\rho_b$ at the zero point.


\bigskip\noindent
\textbf{2.6.}
In this subsection we will derive Theorem \ref{2.5} from 
Proposition \ref{2.45}. Let us consider the factor of product \eqref{baddit}
corresponding to any pair of indices $a<b\ts$.
There $k=1\lcd d\ts$.
Since the indices $i_1\lcd i_d$ and $j_1\lcd j_d$ are pairwise distinct, 
we may assume that $d\le\nu_a$ and $d\le\nu_b$ in \eqref{baddit}.
If $\nu_a<\nu_b$ then $k<\nu_b\ts$.
Therefore any factor of \eqref{baddit}
corresponding to $a<b$ with $\nu_a<\nu_b$ is regular
at any point $\mu\in\t_m^*$ where the weight $\la+\rho$ is dominant.

If $\nu_a\ge\nu_b$ then the denominator in \eqref{baddit}
is zero if and only if $k=d=\nu_b$ and 
$
\la_a+\rho_a=\la_b+\rho_b\,.
$
Hence any factor in the product \eqref{baddit} with 
$\nu_a\ge\nu_b$ has a simple pole as a function
of $\la_a-\la_b+\rho_a-\rho_b\,$, at the zero~point.
Moreover, the residue of the latter function at the zero point equals
$$
\frac1{(\nu_b-1)\,!}\ \,
\sum_{i_1,\ldots,i_{\nu_b}}\ 
P_{\ts i_11}\ts\ldots\ts P_{\ts i_{\nu_b}\nu_b}
$$
where the indices $i_1\lcd i_{\nu_b}$ are distinct 
and taken from 
$1\lcd\nu_a\ts$.

Using the observation on \eqref{beigen}
made at the end of the previous subsection,
we now complete the proof of Theorem \ref{2.5}
for any dominant weight $\la+\rho\ts$.

Note that Proposition~\ref{2.4}, like Theorem \ref{2.5} above,
can also be derived from Proposition \ref{2.45}\tts; see 
the proof of \cite[Proposition 4.4]{KN}. Further,
Proposition~\ref{2.45} and Theorem \ref{2.5} both
can be proved by using
\cite[Subsections 1.4 and~4.4]{KN}.   


\bigskip\noindent
\textbf{2.7.}
In this subsection we will give another proof of Theorem \ref{2.5}.
It provides a multiplicative formula for the operator $I(\mu)$
with dominant $\la+\rho$. Consider
the product \eqref{bmu}. It is taken over the pairs
$(p\com q)$ where $p<q$ while $p$~and~$q$ belong to
different segments of the sequence $1\lcd N\ts$; see \eqref{pnu}. 
Let $a$ and $b$ be the numbers of these two segments, so that
$a<b\ts$. Let us now rearrange the pairs $(p\com q)$ in the product 
\eqref{bmu} as follows.
 
The new order on the pairs $(p\com q)$
will be an extension of the lexicographical order
on the corresponding pairs $(a\com b)\ts$.
To define the extension, we have to order the pairs $(p\com q)$
corresponding to a given $(a\com b)\ts$. 
Take another pair $(r\com s)$ such that
the indices $r$ and $s$ belong to the segments $a$ and $b\ts$,
that is to the same segments as the indices $p$ and $q$ respectively.
For $\nu_a<\nu_b\ts$,
the pair $(p\com q)$ will precede $(r\ts,s)$ if $p<r$
or if $p=r$ and $q>s\ts$.
For $\nu_a\ge\nu_b\ts$,
the pair $(p\com q)$ will precede $(r\ts,s)$ if $q>s$
or if $q=s$ and $p<r\ts$. 
By exchanging commuting factors in \eqref{bmu}, 
this rearrangement does not alter the value of the ordered~product.

Let $i$ and $j$ be the numbers of the indices $p$ and $q$
within their segments, so that by definition 
we have the equalities \eqref{xp} and
\begin{equation}
\label{xq}
\textstyle
x_q=\mu_b+\rho_b+\frac12+\nu_b-j
\ts.
\end{equation}
Consider the factor $R_{pq}(x_q-x_p)$ in the reordered product
\eqref{bmu}. As a function of $\mu\ts$, this factor
has a pole at $x_p=x_q\ts$. The latter equation can be written as
\begin{equation}
\label{lala}
\la_a-\la_b+\rho_a-\rho_b=i-j
\end{equation}
using \eqref{xp} and \eqref{xq}. Hence if $x_p=x_q$ while
$\la+\rho$ is dominant, then $i\ge j\ts$.

\enlargethispage{2pt}

First, suppose that $\nu_a<\nu_b\ts$. If $i\ge j\ts$, then $j<\nu_b$
since $i\le \nu_a\ts$. Then the index $q$ is not the last
in its segment, so that the pair $(p\com q+1)$ immediately 
precedes $(p\com q)$ in the new ordering.
Consider the pairs which follow
$(p\com q)$ in the new ordering. Take the product of the
factors in \eqref{bmu} corresponding to the latter pairs,
and multiply it on the right by
\begin{equation}
\label{qq}
R_{q,q+1}\ts(x_q-x_{q+1})=R_{q,q+1}\ts(1)=1-P_{q,q+1}\,.
\end{equation}
Using \eqref{ybe}, the resulting product
is also divisible by \eqref{qq} on the left.
Due to \eqref{nopole} 
we can therefore replace in \eqref{bmu} 
the product of two adjacent factors 
$$
R_{p,q+1}\ts(x_{q+1}-x_p)\,
R_{pq}\ts(x_q-x_p)
\quad\ \text{by}\ \quad
1-(\ts P_{pq}+P_{p,q+1}\ts)/(x_{q+1}-x_p)
$$
without changing the restriction of the operator
\eqref{bmu} to the subspace   
\eqref{subspace}. But the replacement
does not have a pole at $x_p=x_q\ts$.
So the factors in \eqref{bmu}
corresponding to the pairs $(p\com q)$ with $\nu_a<\nu_b$
do not increase the order of the
pole of $I(\mu)$ at any point $\mu$
such that the weight $\la+\rho$ is dominant.

Next, suppose that $\nu_a\ge\nu_b$ while $i>1\ts$,
so that the index $p$ is not the first in its segment. 
Then the pair $(p-1\com q)$ immediately precedes
$(p\com q)$ in the new ordering. Consider the pairs following
$(p\com q)\ts$. Take the product of the
factors in \eqref{bmu} corresponding to the latter pairs.
Multiply it on the right~by
\begin{equation}
\label{pp}
R_{p-1,p}\ts(x_{p-1}-x_p)=R_{p-1,p}\ts(1)=1-P_{p-1,p}\,.
\end{equation}
Using \eqref{ybe}, the resulting product
is also divisible by \eqref{pp} on the left.
Due to \eqref{nopole} 
we can now replace in \eqref{bmu} 
the product of two adjacent factors 
$$
R_{p-1,q}\ts(x_q-x_{p-1})\,
R_{pq}\ts(x_q-x_p)\,
\quad\ \text{by}\ \quad
1-(\ts P_{p-1,q}+P_{pq}\ts)/(x_q-x_{p-1})
$$
without changing the restriction of \eqref{bmu} to 
\eqref{subspace}. The replacement has no pole at $x_p=x_q\ts$.
So the factors in \eqref{bmu}
corresponding to the pairs $(p\com q)$ with $\nu_a\ge\nu_b$ and $i>1$
do not increase the order of the
pole of $I(\mu)$ at any point $\mu\ts$.

Last, suppose that $\nu_a\ge\nu_b$ and $i=1\ts$.
If $x_p=x_q$ then $\la_a+\rho_a=\la_b+\rho_b$
and $j=1$ whenever the weight $\la+\rho$ is dominant, due to \eqref{lala}.
The observation on (2.24) made at the end of Subsection 2.5
now proves Theorem \ref{2.5}. 

Note that all the replacements in the product
\eqref{bmu} described above can be made simultaneously.
Hence our argument provides an explicit 
formula for the operator $I(\mu)$ whenever $\la+\rho$ is dominant.
See also \cite[Subsection~2.3]{GNP}.


\bigskip\noindent
\textbf{2.8.}
In this subsection we will generalize Theorem \ref{1.1}.
Theorem \ref{2.5} allows us
to determine the intertwining operator $I(\mu)$ 
of the $\Y(\gl_n)\ts$-modules \eqref{intoper1}
for any $\mu\in\t_m^*\ts$, provided the weight $\la+\rho$ is dominant.
Our generalization of Theorem \ref{1.1} is based on the following
lemma. For any index $c=1\lcd m-1$ let $s_c\in S_m$ 
be the transposition of $c$ and $c+1\ts$. Here
the symmetric group $S_m$ acts on the numbers $1\lcd m$ by
their permutations. The latter correspond to
permutations of the basis vectors $E_{11}\lcd E_{mm}$ of $\t_m\ts$.

\begin{lemma*}
\label{2.55}
Fix\/ $c>0$ and suppose that both\/ $\la+\rho$ and\/ $s_c(\la+\rho)$
are dominant. Then the images of 
the intertwining operators $I(\mu)$ and $I(s_c\circ\mu)$ 
corresponding to the pairs $(\la\com\mu)$ and 
$(s_c\circ\la\com s_c\circ\mu)$ 
are equivalent as\/ $\Y(\gl_n)\ts$-modules.
\end{lemma*}

\begin{proof}
Note that in this lemma the intertwining operator $I(s_c\circ\mu)$ 
corresponds to the weight $s_c\circ\la\ts$, not to $\la\ts$.
Moreover, the source and target $\Y(\gl_n)\ts$-modules
of this operator are different from those of the operator
$I(\mu)$ in general. This difference
should not cause any confusion however.

Let $\lac$ and $\muc$ and $\rhoc$ be the weights of $\gl_2$ with
the labels $\la_c\com\la_{c+1}$ and $\mu_c\com\mu_{c+1}$ and 
$\rho_c\com\rho_{c+1}$ respectively. The weights
$\lac+\rhoc$ and $s_1(\lac+\rhoc)$ of $\gl_2$ are dominant. By
using Theorem \ref{2.5} with $m=2$ we get the 
$\Y(\gl_n)\ts$-intertwining operators
\begin{align}
\label{fint}
\Ph^{\,\nu_c}_{\mu_c+\rho_c+\frac12}
\ot
\Ph^{\,\nu_{c+1}}_{\mu_{c+1}+\rho_{c+1}+\frac12}
&\to\,
\Ph^{\,\nu_{c+1}}_{\mu_{c+1}+\rho_{c+1}+\frac12}
\ot
\Ph^{\,\nu_c}_{\mu_c+\rho_c+\frac12}\,,
\\[8pt]
\label{gint}
\Ph^{\,\nu_{c+1}}_{\mu_{c+1}+\rho_{c+1}+\frac12}
\ot
\Ph^{\,\nu_c}_{\mu_c+\rho_c+\frac12}
&\to\,
\Ph^{\,\nu_c}_{\mu_c+\rho_c+\frac12}
\ot
\Ph^{\,\nu_{c+1}}_{\mu_{c+1}+\rho_{c+1}+\frac12}\,.
\end{align}
These two operators are inverse to each other.
This assertion can be proved first for any generic weight 
$\mu$ of $\gl_m\ts$,
either by a direct calculation employing the definition \eqref{bmu}, or
by observing that for any generic $\mu$
all the double tensor products above are $\Y(\gl_n)\ts$-irreducible
while 
\eqref{fint} and \eqref{gint} respectively map
$$
\ph_{\nu_c}\ot\ph_{\nu_{c+1}}
\mapsto\ts
\ph_{\nu_{c+1}}\ot\ph_{\nu_c}
\quad\ \text{and}\ \quad
\ph_{\nu_{c+1}}\ot\ph_{\nu_c}
\mapsto\ts
\ph_{\nu_c}\ot\ph_{\nu_{c+1}}\,.
$$
Then the assertion extends to $\mu$
such that $\la+\rho$ and $s_c(\la+\rho)$ are dominant.

Now denote by $I$ the operator which acts 
on the tensor product of $c\,$th and
$(c+1)\,$th factors of \eqref{basprod1}
as the intertwining operator \eqref{fint}, and
which acts trivially on all other $m-2$ tensor factors of \eqref{basprod1}.
Arguing like above, that is either performing a direct calculation,
or using the irreducibility of the source and target
$\Y(\gl_n)\ts$-modules in \eqref{intoper1} for any generic $\mu$,
we get the relation 
$$
I(\mu)=I\,I(s_c\circ\mu)\,I\,. 
$$
It proves the lemma, since $I$ is invertible
and intertwines $\Y(\gl_n)\ts$-modules.
\qed
\end{proof}

For any $\la\in\t_m^*$ denote by $S_\la$
the subgroup of $S_m$ consisting of all elements $w$ such that 
$w\circ\la=\la\ts$. Let $\mathcal{O}$ be any orbit of the shifted 
action of the subgroup $S_\la\subset S_m$ on $\t_m^*\ts$.
If $\nu_1\lcd\nu_m\in\{\tts1\lcd n-1\}$ for at least one
weight $\mu\in\mathcal{O}$,
then every $\mu\in\mathcal{O}$
satisfies the same condition. 
Suppose this is the case for $\mathcal{O}$.
If $\la+\rho$ is dominant,
then there is at least one weight $\mu\in\mathcal{O}$
such that the pair $(\la\com\mu)$ is good.
Theorem \ref{1.1} generalizes due to the following proposition.

\begin{proposition*}
\label{2.575}
If $\la+\rho$ is dominant, then for all 
$\mu\in\mathcal{O}$
the images of the corresponding operators\/ $I(\mu)$
are equivalent to each other as\/ $\Y(\gl_n)\ts$-modules.
\end{proposition*}

\begin{proof}
Take any $w\in S_\la\ts$ and any reduced decomposition 
$w=s_{c_l}\ldots s_{c_1}\ts$. It can be derived from
\cite[Corollary VI.1.2]{B} that
the weight  $s_{c_k}\ldots s_{c_1}(\la+\rho)$ of $\gl_m$ is
dominant for each $k=1\lcd l\ts$.
Proposition \ref{2.575} now follows by applying Lemma \ref{2.55} repeatedly.
Note that this proposition can also be proved by
using the results of Zelevinsky \cite[Theorem 6.1]{Z}
together with those of \cite{AS,D1}.
\qed
\end{proof}

Thus all assertions of Theorem \ref{1.1} will remain valid
if we replace the good pair there by any pair
$(\la\com\mu)$ such that the weight $\la+\rho$ of $\gl_m$ is dominant.   
However, we still assume that $\nu_1\lcd\nu_m\in\{\tts1\lcd n-1\}$
for the latter pair.


\section*{\normalsize\bf 3. More intertwining operators}
\label{S3}
\setcounter{section}{3}
\setcounter{theorem*}{0}
\setcounter{equation}{0}


\textbf{3.1.}
Now let $\la$ and $\mu$ be any weights of the Lie algebra
$\f_m=\sp_{2m}$ or $\f_m=\so_{2m}$
such that all labels of the weight 
$\nu=\la-\mu+\ka$ are in $\{\tts1\lcd n-1\}\ts$. We will
keep using the notation \eqref{pnu},\eqref{xp},\eqref{bmu}. 
But now $\la_a\com\mu_a\com\nu_a$ and $\rho_a$ with $a=1\lcd m$
are labels of weights of $\f_m\ts$. Recall that 
$\ka_a=n/2$ by definition.  

Let the weight $\mu\in\h_m^*$ vary while $\nu$ is fixed.  
Determine the rational function $B(\mu)$ by the same
formula \eqref{bmu} as for $\Y(\gl_n)\ts$.
Also take the ordered product
\begin{equation}
\label{cmu}
C(\mu)=\prod_{(p,q)}^{\longleftarrow}
\,\Rt_{\ts pq}\ts(\tts n-x_p-x_q\tts)
\end{equation}
where $1\le p<q\le N$ and the pairs $(p\com q)$ are ordered 
lexicographically. Here the reversed arrow indicates that the factors
corresponding to these pairs are arranged from right to left.
Note that $C(\mu)$ is a rational function of $\mu\in\h_m^*$
without poles at the generic weights $\mu$ of $\f_m\ts$.

Take the sequence $1''\lcd N\tts''$ 
introduced in the previous subsection.
Let $Q_\nu$ be the 
linear operator on $(\CC^n)^{\tts\ot\tts N}$ which
for each $p=1\lcd N$ exchanges the tensor factors $\CC^n$
labeled by $p$ and $p\ts''\ts$.
Then $Q_\nu=P$ if $m=1$ and $\nu_1=2\ts$.

\begin{proposition*}
\label{2.6}
Suppose the weight $\mu$ of\/ $\f_m$ is generic.
Then\/ $Q_\nu\ts B(\mu)\,C(\mu)$ is an intertwining operator of\/ 
$\Y(\g_n)\ts$-modules
\begin{equation}
\label{bc}
\Ph^{\tts1}_{\tts x_1}
\ot\ldots\ts\ot\ts
\Ph^{\tts1}_{\tts x_N}
\to\,
\Ph^{\ts-1}_{\tts x_{1''}}
\ot\ldots\ts\ot\ts
\Ph^{\ts-1}_{\tts x_{N''}}
\ts.
\end{equation}
\end{proposition*}

\begin{proof}
Under the action of $\Y(\gl_n)$ on the source module in \eqref{bc},
\begin{equation*}
S(x)
\,\mapsto\,
\Rh_{\ts0N}(\ts x_N+x)
\,\ldots\,\ts
\Rh_{\ts01}(x_1+x)
\,\,
\Rb_{\ts01}(x_1-x)
\,\ldots\,\ts
\Rb_{\ts0N}(\ts x_N-x)\ts.
\end{equation*}
Let us denote by $Y$ the right hand side of 
this assignment.
Here we use the subscripts $0\com1\lcd N$ rather than
$1\com2\lcd N+1$ to label the tensor factors of
$(\ts\End\CC^n)^{\tts\ot\tts(N+1)}\ts$,
like we did in the proof of Lemma \ref{2.1}. Further, let us denote
$$
Y'=
\Rb_{\ts01}(x_1-x)
\,\ldots\,\ts
\Rb_{\ts0N}(\ts x_N-x)
\,\ts
\Rh_{\ts0N}(\ts x_N+x)
\,\ldots\,\ts
\Rh_{\ts01}(x_1+x)\ts.
$$
By using the relation \eqref{bybe} repeatedly, we get 
the equality $C(\mu)\,Y=Y'\ts C(\mu)\ts$.
Here we also use the relation $P\ts\Rt(x)\ts P=\Rt(x)\ts$.

Under the action of $\Y(\g_n)$ on the target module in \eqref{bc},
\begin{equation*}
S(x)
\,\mapsto\,
\Rb_{\ts0N}(\ts x_{N''}-x)
\,\ldots\,\ts
\Rb_{\ts01}(x_{1''}-x)
\,
\Rh_{\ts01}(x_{1''}+x)
\,\ldots\,\ts
\Rh_{\ts0N}(\ts x_{N''}+x)\ts.
\end{equation*}
Denote by $Y''$ the right hand side of the latter assignment.
Observe that $k\ts''=(N-k+1)\ts'$
for each $k=1\lcd N\ts$. Therefore we can write
$$
Y''=
\Rb_{\ts0N}(\ts x_{1'}-x)
\,\ldots\,\ts
\Rh_{\ts01}(x_{N'}-x)
\,
\Rh_{\ts01}(x_{N'}+x)
\,\ldots\,\ts
\Rh_{\ts0N}(\ts x_{1'}+x)\ts.
$$
But by using the relations \eqref{tybe} and \eqref{hybe} repeatedly, we get
\begin{gather*}
Q_\nu\ts B(\mu)\,Y'=
\\[2pt]
Q_\nu\,
\Rb_{\ts01'}(x_{1'}-x)
\,\ldots\,\ts
\Rb_{\ts0N'}(\ts x_{N'}-x)
\,
\Rh_{\ts0N'}(\ts x_{N'}+x)
\,\ldots\,\ts
\Rh_{\ts01'}(x_{1'}+x)\,B(\mu)
\\[6pt]
=\,Y''\ts Q_\nu\ts B(\mu)\,;
\end{gather*}
see also the proof of Proposition \ref{2.3}.
Thus we get the equality
$$
Q_\nu\ts B(\mu)\,C(\mu)\,Y=\ts Y''\ts Q_\nu\ts B(\mu)\,C(\mu)\ts,
$$
which proves the claim.
\qed
\end{proof}


\noindent
\textbf{3.2.}
Observe that $\Pt\,P=\pm\ts\Pt\ts$.
Thus $\Pt\,R(1)=0$ in the case $\g_n=\so_n\ts$.
Hence in this case
the restriction of the operator $\Pt_{pq}$ to the subspace
\eqref{subspace} is zero for 
any two distinct indices $p\com q$ from the same segment of the sequence
$1\lcd N\ts$. Here we mean the segments of lengths $\nu_1\lcd\nu_m$
as defined in Subsection 2.3.

In the case $\g_n=\so_n$ the relation \eqref{cybe} now implies 
that when considering the restriction of
the operator $C(\mu)$ to the subspace \eqref{subspace},
we can skip those factors in the product \eqref{cmu} 
which correspond to the pairs $(p,q)$ where both $p$ and $q$
belong to the same segment: \ts skipping does not 
change the restriction.  
In particular, in the case $\g_n=\so_n$ the restriction of
the operator $C(\mu)$ to the subspace \eqref{subspace}
does not have a pole 
if $\mu_a-\mu_b\notin\ZZ$ and $\mu_a+\mu_b\notin\ZZ$
whenever $a\neq b\ts$;
the condition $2\ts\mu_a\notin\ZZ$ is not needed here.

Let us give an analogue of Proposition \ref{2.4}
for $C(\mu)\ts$. 
Except in the proof of Lemma \ref{1.15}, we worked with
any symmetric or alternating non-degenerate 
bilinear form $\langle\ ,\,\rangle$ on $\CC^n$ so far.
Choose the form as in the proof, so that
\eqref{teij} holds. The elements
$
E_{\ts ij}-\widetilde{E}_{\ts ij}
$
with $i\le j$ span a Borel subalgebra of $\g_n\subset\gl_n$ then\tts; 
the elements
$
E_{\ts ii}-\widetilde{E}_{\ts ii}
$
span~the corresponding Cartan subalgebra~of~$\g_n\ts$.

\begin{proposition*}
\label{2.7}
For any generic weight $\mu$ of\/ $\f_m$ the vector\/
$\ph_{\nu_1}\ns\ot\ts\ldots\ts\ot\ts\ph_{\nu_m}\ns$
is an eigenvector of the operator $C(\mu)$ on
$(\CC^n)^{\tts\ot\tts N}$. The eigenvalue is the product
\begin{equation}
\label{ceigen}
\prod_{1\le a<b\le m}\ 
\left\{
\begin{array}{cl}
\displaystyle
\frac{\,\la_a+\la_b+\rho_a+\rho_b}{\,\mu_a+\mu_b+\rho_a+\rho_b}
&\quad\textrm{if}\quad\ \nu_a+\nu_b\ge n\,;
\\[16pt]
\displaystyle
1
&\quad\textrm{if}\quad\ \nu_a+\nu_b\le n\,;
\end{array}
\right.   
\end{equation}
multiplied in the case of\/ $\f_m=\sp_{2m}$ by the product
\begin{equation}
\label{deigen}
\prod_{1\le a\le m}\ 
\left\{
\begin{array}{cl}
\displaystyle
\frac{\,\la_a+\rho_a}{\,\mu_a+\rho_a}
&\quad\textrm{if}\quad\ 2\ts\nu_a\ge n\,;
\\[16pt]
\displaystyle
1
&\quad\textrm{if}\quad\ 2\ts\nu_a\le n\,.
\end{array}
\right.   
\end{equation}
\end{proposition*}

\begin{proof}
This proposition immediately follows from its particular cases
of $m=1$ and $m=2\ts$. We will consider these two cases only. 
First suppose that $m=1\ts$. In this case, for $\g_n=\so_n$
the restriction of the operator $C(\mu)$ to the subspace \eqref{subspace}
is the identity operator\ts; see the observation made in the very
beginning of the present subsection. Suppose that $\g_n=\sp_n\ts$.
Put $\ep=(-1)^{\ts[\ts\nu_1/2\ts]}\ts$. Then
\begin{gather}
\nonumber
C(\mu)\,\ph_{\nu_1}
\,=\!
\prod_{i=1\lcd\nu_1-1}^{\longleftarrow}
\!\Bigl(\,
\prod_{j=i+1,\ldots,\nu_1}^{\longleftarrow}
\Rt_{ij}\ts(\tts n-x_i-x_j\tts)
\Bigr)\ \times
\\[6pt]
\nonumber
\ep\,A_{\nu_1}
(\tts e_{\nu_1}\ot\ts\ldots\ts\ot e_1\tts)\,=\,
\\[6pt]
\label{disvec}
\ep\,A_{\nu_1}
\prod_{i=1\lcd\nu_1-1}^{\longrightarrow}
\Bigl(\,\,
\prod_{j=i+1,\ldots,\nu_1}^{\longrightarrow}
\Rt_{ij}\ts(\tts n-x_i-x_j\tts)
\Bigr)\,
(\tts e_{\nu_1}\ot\ts\ldots\ts\ot e_1\tts)
\end{gather}
where the latter equality is obtained by using \eqref{cybe} and \eqref{ak}.
We will prove by induction on $\nu_1$ that the vector \eqref{disvec}
equals $\ph_{\nu_1}$ 
multiplied by the scalar~\eqref{deigen} where $m=1\ts$.
Recall that $\nu_1=\la_1-\mu_1+n/2\ts$.
If $\nu_1\le n/2$ then here we have
$$
\Pt_{ij}(\tts e_{\nu_1}\ot\ts\ldots\ts\ot e_1\tts)=0
$$
by our choice of the form $\langle\ ,\,\rangle$ on $\CC^n$.
Hence then the vector \eqref{disvec}
equals $\ph_{\nu_1}$ as required.
We will also use that equality for $\nu_1\le n/2$ as the induction base.

Let $\nu_1>n/2\ts$. In particular, then $\nu_1>1$ because $n$ is even.
The induction assumption then implies that
\begin{gather*}
(1\ot A_{\ts\nu_1-1})
\prod_{i=2,\ldots,\nu_1-1}^{\longrightarrow}
\Bigl(\,\,
\prod_{j=i+1,\ldots,\nu_1}^{\longrightarrow}
\Rt_{ij}\ts(\tts n-x_i-x_j\tts)
\Bigr)\,
\\[6pt]
(\tts e_{\nu_1}\ot e_{\nu_1-1}\ot\ts\ldots\ts\ot e_1\tts)\,=\,
u\,e_{\nu_1}\ot A_{\ts\nu_1-1}\ts
(\tts e_{\nu_1-1}\ot\ts\ldots\ts\ot e_1\tts)
\end{gather*}
where
$$
u\,=\,
\frac{\,\la_1+\rho_1-1}{\,\mu_1+\rho_1}\,\,.
$$
We use the inequality $2\ts(\nu_1-1)\ge n$ which follows 
from $\nu_1>n/2\ts$, because $n$ is even.
Arguing like in the proof of Lemma \ref{2.1} and using
the relation \eqref{nopolet}~we~get 
\begin{gather*}
(\tts1\ot A_{\ts\nu_1-1}\tts)
\prod_{j=2,\ldots,\nu_1}^{\longrightarrow}
\!\Rt_{\ts1j}\ts(\tts n-x_1-x_j\tts)\,=\,
\\
(\tts1\ot A_{\ts\nu_1-1}\tts)
\Bigl(
\ts1-\ts
(\tts n-x_1-x_2\tts)^{-1}\,
\sum_{j=2}^{\nu_1}\,
\Pt_{\ts 1j}\Bigr)\,.
\end{gather*}
Further,
by our choice of the form $\langle\ ,\,\rangle$ on $\CC^n$ the vector
$$
A_{\nu_1}\Pt_{\tts1j}\ts(\tts e_{\nu_1}\ot\ts\ldots\ts\ot e_1\tts)
$$
is equal to
$$
2\ts A_{\nu_1}(\tts e_{\nu_1}\ot\ts\ldots\ts\ot e_1\tts)
$$
if $j=2\ts\nu_1-n\ts$, and is equal to zero for any other $j\ge 2\ts$.
By writing \eqref{disvec} as
\begin{gather*}
\ep\,A_{\nu_1}
\prod_{j=2,\ldots,\nu_1}^{\longrightarrow}
\!\Rt_{\ts1j}\ts(\tts n-x_1-x_j\tts)\ \times
\\[2pt]
\prod_{i=2\lcd\nu_1-1}^{\longrightarrow}
\Bigl(\,\,
\prod_{j=i+1,\ldots,\nu_1}^{\longrightarrow}
\Rt_{ij}\ts(\tts n-x_i-x_j\tts)
\Bigr)\ 
(\tts e_{\nu_1}\ot\ts\ldots\ts\ot e_1\tts)
\end{gather*}
we now see that the vector
\eqref{disvec} is equal to $\,\ph_{\nu_1}$ multiplied by
the scalar
$$
\frac{n-x_1-x_2-2}{n-x_1-x_2}\,\,u\,=\,
\frac{\,\la_1+\rho_1}{\,\mu_1+\rho_1}\,\,.
$$

Thus we have finished the proof of Proposition \ref{2.7}
in the case $m=1\ts$, and will now suppose that $m=2\ts$. 
Then by the definition \eqref{cmu} the operator $C(\mu)$
is the ordered product of $\Rt_{\ts pq}\ts(\tts n-x_p-x_q\tts)$
over the pairs $(p\com q)$ where $1\le p<q\le\nu_1+\nu_2\ts$;
the arrangement of the pairs is reversed lexicographical.
Without changing the product,
we can rearrange these pairs as follows.
From left to right, first come the pairs $(p\com q)$ where
$\nu_1<p<q\le\nu_1+\nu_2\ts$, second come the pairs
where $1\le p\le \nu_1$ and $\nu_1<q\le\nu_1+\nu_2\ts$;
third come the pairs where $1\le p<q\le\nu_1\ts$.
Within each of the three groups,
the arrangement of the pairs $(p\com q)$ is still reversed lexicographical.

Consider the two products of $\ts\Rt_{\ts pq}\ts(\tts n-x_p-x_q\tts)\ts$,
over the first and the third group of pairs
$(p\com q)\ts$. The already settled case $m=1$ implies
that $\ph_{\nu_1}\ot\ph_{\nu_2}$ is an eigenvector
for each of these two products. If $\g_n=\so_n$ then
each of the two corresponding eigenvalues is $1\ts$.
If $\g_n=\sp_n$ then the product of the two eigenvalues
is \eqref{deigen} where $m=2\ts$.
For $\g_n=\so_n\com\sp_n$
we will show that $\ph_{\nu_1}\ot\ph_{\nu_2}$
is an eigenvector of the product of 
$\Rt_{\ts pq}\ts(\tts n-x_p-x_q\tts)$
over the second group of pairs,
with the eigenvalue \eqref{ceigen} where $m=2\ts$.
This will settle the case $m=2\ts$.

Use the induction on $\nu_1\ts$. 
Denote the last mentioned product by $Z$, so that
\begin{equation}
\label{asreq}
Z\ =\prod_{i=1,\ldots,\nu_1}^{\longleftarrow}
\Bigl(\,\,
\prod_{j=1,\ldots,\nu_2}^{\longleftarrow}
\Rt_{\ts i,\nu_1+j}\ts(\tts n-x_i-x_{\nu_1+j}\tts)
\Bigr)
\end{equation}
by definition. By using \eqref{cybe} and \eqref{ak}, 
we get the equality
\begin{gather}
\nonumber
Z\,(A_{\nu_1}\ot A_{\nu_2})\,
(\tts e_{\nu_1}\ot\ts\ldots\ts\ot e_1
\ot 
e_{\nu_2}\ot\ts\ldots\ts\ot e_1\tts)
\,=\,
\\[12pt]
(A_{\nu_1}\ot A_{\nu_2})
\label{anu}
\prod_{i=1,\ldots,\nu_1}^{\longrightarrow}
\Bigl(\,\,
\prod_{j=1,\ldots,\nu_2}^{\longrightarrow}
\Rt_{\ts i,\nu_1+j}\ts(\tts n-x_i-x_{\nu_1+j}\tts)
\Bigr)
\\[8pt]
\nonumber
(\tts e_{\nu_1}\ot\ts\ldots\ts\ot e_1
\ot 
e_{\nu_2}\ot\ts\ldots\ts\ot e_1\tts)\,.
\end{gather}
If $\nu_1+\nu_2\le n$ then
by our choice of the form $\langle\ ,\,\rangle$ on $\CC^n$
$$
\Pt_{\ts i,\nu_1+j}\ts
(\tts e_{\nu_1}\ot\ts\ldots\ts\ot e_1
\ot 
e_{\nu_2}\ot\ts\ldots\ts\ot e_1\tts)=0
$$ 
for $i\com j$ as above. Hence then the product 
\eqref{asreq} acts on the vector 
$\ph_{\nu_1}\ot\ph_{\nu_2}$ as the identity.
We will also use this result for $\nu_1+\nu_2\le n$ as the induction~base.

Let $\nu_1+\nu_2>n\ts$. Then $\nu_1>1$ because $\nu_2<n\ts$. 
By the induction~assumption
\begin{gather*}
(1\ot A_{\ts\nu_1-1}\ot A_{\nu_2})
\prod_{i=2,\ldots,\nu_1}^{\longrightarrow}
\Bigl(\,\,
\prod_{j=1,\ldots,\nu_2}^{\longrightarrow}
\Rt_{\ts i,\nu_1+j}\ts(\tts n-x_i-x_{\nu_1+j}\tts)
\Bigr)\,
\\[6pt]
(\tts e_{\nu_1}\ot e_{\nu_1-1}\ot\ts\ldots\ts\ot e_1
\ot 
e_{\nu_2}\ot\ts\ldots\ts\ot e_1\tts)\,=\,
\\[16pt]
v\,e_{\nu_1}\ot (A_{\ts\nu_1-1}\ot A_{\nu_2})\,
(\tts e_{\nu_1-1}\ot\ts\ldots\ts\ot e_1
\ot 
e_{\nu_2}\ot\ts\ldots\ts\ot e_1\tts)
\end{gather*}
where
\begin{equation*}
v\,=\,
\frac{\,\la_1+\la_2+\rho_1+\rho_2-1}{\,\mu_1+\mu_2+\rho_1+\rho_2}\,\,.
\end{equation*}
\\[-2pt]
Arguing like in the proof of Lemma \ref{2.1} and using
the relation \eqref{nopolet}~we~get
the equality in the algebra $(\ts\End\CC^n)^{\tts\ot\tts(\nu_1+\nu_2)}$ 
\begin{gather*}
(\tts1\ot A_{\ts\nu_2}\tts)
\prod_{j=1,\ldots,\nu_2}^{\longrightarrow}
\!\Rt_{\ts1,\nu_1+j}\ts(\tts n-x_1-x_{\nu_1+j}\tts)\,=\,
\\
(\tts1\ot A_{\ts\nu_2}\tts)\,
\Bigl(
\ts1-\ts
(\tts n-x_1-x_{\nu_1+1}\tts)^{-1}\,
\sum_{j=1}^{\nu_2}\,
\Pt_{\ts 1,\nu_1+j}\Bigr)\,.
\end{gather*}
Further, 
by our choice of the form $\langle\ ,\,\rangle$ on $\CC^n$ the vector
$$
(A_{\nu_1}\ot A_{\nu_2})\,\Pt_{\tts1q}\ts
(\tts e_{\nu_1}\ot\ts\ldots\ts\ot e_1
\ot 
e_{\nu_2}\ot\ts\ldots\ts\ot e_1\tts) 
$$
is equal to
\begin{equation*}
(A_{\nu_1}\ot A_{\nu_2})\,
(\tts e_{\nu_1}\ot\ts\ldots\ts\ot e_1
\ot 
e_{\nu_2}\ot\ts\ldots\ts\ot e_1\tts)
\end{equation*}
\\[-6pt]
if $j=\nu_1+\nu_2-n\ts$, and is equal to zero for any other $j\ts$.
Writing 
\eqref{anu}~as
\begin{align*}
(A_{\nu_1}\ot A_{\nu_2})
&\prod_{j=1,\ldots,\nu_2}^{\longrightarrow}
\Rt_{\ts1,\nu_1+j}\ts(\tts n-x_1-x_{\nu_1+j}\tts)\ \times
\\[2pt]
\prod_{i=2,\ldots,\nu_1}^{\longrightarrow}
\Bigl(\,\,
&\prod_{j=1,\ldots,\nu_2}^{\longrightarrow}
\Rt_{\ts i,\nu_1+j}\ts(\tts n-x_1-x_{\nu_1+j}\tts)
\Bigr)
\end{align*}
we now see that $\ph_{\nu_1}\ot\ph_{\nu_2}$
is an eigenvector for \eqref{asreq} with the eigenvalue
$$
\frac{n-x_1-x_{\nu_1+1}-1}{n-x_1-x_{\nu_1+1}}\,\,v\,=\,
\frac{\,\la_1+\la_2+\rho_1+\rho_2}{\,\mu_1+\mu_2+\rho_1+\rho_2}\,\,.
$$
This observation completes the proof of Proposition \ref{2.7}.
\qed
\end{proof}


\noindent
\textbf{3.3.}
Using the relations \eqref{ybe},\eqref{cybe} and $R(1)=1-P$ 
one shows that for any generic weight $\mu$ of $\f_m$
the operator $C(\mu)$ preserves the subspace \eqref{subspace},
see the proof of Proposition \ref{2.7}.
Moreover, we have another proposition.
It can be used to give another proof of Proposition \ref{2.7},
see the proof of \cite[Proposition~4.6]{KN}.

\begin{proposition*}
\label{2.75}
For any generic weight\/ $\mu$ of\/ $\f_m$ 
the restriction of\/ 
$C(\mu)$ to the subspace \eqref{subspace} coincides with that of the operator
\\[8pt]
\begin{equation*}
\label{caddit}
\prod_{1\le a\le b\le m}^{\longleftarrow}\ 
\left\{
\begin{array}{l}
\, 
\\[-20pt]
\displaystyle
\,1\ +\   
\sum_{d>0}\ 
\sum_{\substack{i_1,\ldots,i_d\\j_1,\ldots,j_d}}\ \, 
\prod_{k=1}^d\ \ \, 
\frac{\Pt_{\ts i_kj_k}}{\la_a+\la_b+\rho_a+\rho_b-k}
\,\quad\textrm{if}\quad\ a<b\,;
\\[16pt]
\displaystyle
\,1\ +\ 
\sum_{d>0}\ 
\sum_{\substack{i_1,\ldots,i_d\\j_1,\ldots,j_d}}\ \, 
\prod_{k=1}^d\ \ \, 
\frac{\Pt_{\ts i_kj_k}}{2\ts(\ts\la_a+\rho_a-k\ts)}
\ \quad\textrm{if}\quad\ a=b\,,\f_m=\sp_{2m}\ts;
\\[28pt]
\displaystyle
\,1
\,\quad\textrm{if}\quad\ a=b\,,\f_m=\so_{2m}\ts;
\end{array}
\right.   
\end{equation*}
\\[4pt]
here the pairs $(a\com b)$ are ordered lexicographically.
If $a<b$ then 
$i_1\lcd i_d$ \/and\/ $\,j_1\lcd j_d$ are
distinct numbers from the $a\,$th and $b\,$th 
segments of the sequence $1\lcd N$ respectively, taken so that
different are all the sets \eqref{dp}.
If $a=b$ then 
$i_1\com j_1\lcd i_d\com j_d$ are pairwise 
distinct numbers from the $a\,$th
segment of the sequence $1\lcd N$ taken so that
different are all the sets of $d$ unordered pairs
\begin{equation}
\label{dup}
\{i_1\com j_1\}\lcd\{i_d\com j_d\}\ts.
\end{equation}
\end{proposition*}

\begin{proof}
This proposition immediately follows from its particular cases
of $m=1$ and $m=2\ts$. We will consider these two cases only. 
First suppose that $m=1\ts$. We have already observed that
then for $\g_n=\so_n$
the restriction of the operator $C(\mu)$ to the subspace \eqref{subspace}
is the identity operator. Suppose $\g_n=\sp_n\ts$. Then
in the second displayed line in Proposition \ref{2.75} we have
the operator on $(\CC^n)^{\tts\ot\tts\nu_1}$
\begin{equation}
\label{aone}
\sum_{d\ge0}\ 
\sum_{\substack{i_1,\ldots,i_d\\j_1,\ldots,j_d}}\ \, 
\prod_{k=1}^d\ \ \, 
\frac{\Pt_{\ts i_kj_k}}{2\ts x_k-n-1}
\end{equation}
where $i_1\com j_1\lcd i_d\com j_d$ are pairwise 
distinct and taken from $1\lcd\nu_1$
so that different are all the sets of $d$ 
unordered pairs~\eqref{dup}. 
Here we use the equality
\begin{equation}
\label{xk}
x_k=\la_1+\rho_1+(n+1)/2-k
\end{equation}
for any $k\le\nu_1\ts$.
We also assume that $1$ is the only term in \eqref{aone} with $d=0\ts$.
On the other hand, for $m=1$ we can write
\begin{equation}
\label{cone}
C(\mu)\,=\!
\prod_{j=2,\ldots,\nu_1}^{\longleftarrow}
\!\Bigl(\ 
\prod_{i=1,\ldots,j-1}^{\longleftarrow}
\Rt_{ij}\ts(\tts n-x_i-x_j\tts)
\Bigr)\,.
\end{equation} 
Let us now relate two operators on the vector space 
$(\CC^n)^{\tts\ot\tts\nu_1}$
by the symbol~$\,\equiv$ if their actions coincide on
the subspace
\begin{equation}
\label{lanone}
\La^{\nu_1}(\CC^n)
\subset
(\CC^n)^{\tts\ot\tts\nu_1}.
\end{equation}
We will establish the relation
$\,\equiv\,$ between \eqref{aone} and \eqref{cone}
by induction on~$\nu_1\ts$. 

If $\nu_1>1$ then we assume the latter relation holds for $\nu_1-1$ 
instead of $\nu_1\ts$; if $\nu_1=1$ then we are not making any assumption.
Arguing like in the proof of Lemma \ref{2.1}
and then using the induction assumption, \eqref{cone}
is related by~$\,\equiv\,$~to
\begin{gather*}
\Bigl(\ 1\,+\,
\sum_{i=1}^{\nu_1-1}\,
\frac{\Pt_{\ts i\ts\nu_1}}
{x_1+x_{\nu_1}-n}
\ \Bigr)\,
\times
\prod_{j=2,\ldots,\nu_1-1}^{\longleftarrow}
\Bigl(\ 
\prod_{i=1,\ldots,j-1}^{\longleftarrow}
\Rt_{ij}\ts(\tts n-x_i-x_j\tts)
\Bigr)
\\[6pt]
\equiv\ 
\Bigl(\ 1\,+\,
\sum_{i=1}^{\nu_1-1}\,
\frac{\Pt_{\ts i\ts\nu_1}}
{x_1+x_{\nu_1}-n}
\ \Bigr)\,
\times\, 
\Bigl(\ \  
\sum_{d\ge0}\ 
\sum_{\substack{i_1,\ldots,i_d\\j_1,\ldots,j_d}}\ \, 
\prod_{k=1}^d\ \ \, 
\frac{\Pt_{\ts i_kj_k}}{2\ts x_k-n-1}
\ \Bigr)
\end{gather*}
where $i_1\com j_1\lcd i_d\com j_d$ are
distinct and taken from $1\lcd\nu_1-1$
so that different are all corresponding sets 
\eqref{dup}. 
The right hand side of the last relation equals
\begin{gather} 
\label{rone}
\sum_{d\ge0}\ 
\sum_{\substack{i_1,\ldots,i_d\\j_1,\ldots,j_d}}\ \, 
\prod_{k=1}^d\ \ \ 
\frac{\Pt_{\ts i_kj_k}}{2\ts x_k-n-1}
\ \ \ \,+
\\
\label{rtwo}
\sum_{d>0}\ \ \ 
\sum_{l=1}^d\ \, 
\sum_{\substack{i_1,\ldots,i_d\\j_1,\ldots,j_d}}\ \, 
\frac{\,\Pt_{\ts i_l\nu_1}+\Pt_{\ts j_l\nu_1}}
{x_1+x_{\nu_1}-n}
\ \ \ts
\prod_{k=1}^d\ \ 
\frac{\Pt_{\ts i_kj_k}}{2\ts x_k-n-1}
\ \ \ \,+
\\
\label{rthree}
\sum_{d\ge0}\ \,
\sum_{\substack{i_1,\ldots,i_d\\j_1,\ldots,j_d}}\ \, 
\sum_i\ \ \ 
\frac{\,\Pt_{\ts i\ts\nu_1}}
{x_1+x_{\nu_1}-n}
\ \ \ts
\prod_{k=1}^d\ \ 
\frac{\Pt_{\ts i_kj_k}}{2\ts x_k-n-1}
\end{gather}
where 
$i$ is taken from $1\lcd\nu_1-1$ 
and is different from $i_1\com j_1\lcd i_d\com j_d\,$.

Consider the sum displayed in the line \eqref{rtwo}. Here we have
\begin{gather*}
\Pt_{\ts i_l\nu_1}
\ \prod_{k=1}^d\ \, 
\Pt_{\ts i_kj_k}
\,=\,
\Bigl(\
\prod_{k\neq l}\, 
\Pt_{\ts i_kj_k}
\Bigr)\,
\Pt_{\ts i_l\nu_1}\,
\Pt_{\ts i_lj_l}\,
\,=
\\
\Bigl(\ 
\prod_{k\neq l}\, 
\Pt_{\ts i_kj_k}
\Bigr)\,
\Pt_{\ts i_l\nu_1}\,
P_{\ts j_l\nu_1}\,
\,\equiv\,
-\ \Bigl(\ 
\prod_{k\neq l}\, 
\Pt_{\ts i_kj_k}
\Bigr)\,
\Pt_{\ts i_l\nu_1}\,
\end{gather*}
where 
the right hand side does not involve 
$j_l\ts$. Similarly, in \eqref{rtwo} we have
\begin{gather*}
\Pt_{\ts j_l\nu_1}
\ \prod_{k=1}^d\ \, 
\Pt_{\ts i_kj_k}
\,\equiv\,
-\ \Bigl(\ 
\prod_{k\neq l}\, 
\Pt_{\ts i_kj_k}
\Bigr)\,
\Pt_{\ts j_l\nu_1}\,
\end{gather*}
where 
the right hand side does not
involve the index $i_l\ts$.

Now fix a number $i\in\{1\lcd\nu_1-1\}$
and take any set of $d$ pairs \eqref{dup}
such that one of the pairs contains the number $i\ts$.
Then $i=i_l$ or $i=j_l$ for some $l\ts$.
Let $j$ be the element of
the pair $\{i_l\com j_l\}$ different from $i$, so that
$j=j_l$ or $j=i_l$ respectively.
If the set of the $d-1$ pairs $\{i_k\com j_k\}$ with $k\neq l$
is also fixed, then $j$ ranges over a set
of cardinality $\nu_1-2\ts d\ts$, namely over the fixed set 
$$
\{1\lcd\nu_1-1\}
\setminus
\{i_1\com j_1\lcd i_{l-1}\com j_{l-1}
\com i\com
i_{l+1}\com j_{l+1}\lcd i_d\com j_d\}\ts.
$$

Let us perform the summation over the indices $i_l\com j_l$ and $l$ 
in \eqref{rtwo} first of all the running indices. After that
rename the running indices $i_{l+1}\lcd i_d$ and $j_{l+1}\lcd j_d$ 
respectively by $i_l\lcd i_{d-1}$ and $j_l\lcd j_{d-1}\ts$.
By the arguments given in the previous two paragraphs, the sum \eqref{rtwo}
gets related by $\,\equiv\,$ to 
\begin{equation}
\label{rfour}
\sum_{d>0}\,
\sum_{\substack{i_1,\ldots,i_{d-1}\\j_1,\ldots,j_{d-1}}}\, 
\sum_i\ \
\frac{2\ts d-\nu_1}
{2\ts x_d-n-1}
\ \ 
\frac{\,\Pt_{\ts i\ts\nu_1}}
{x_1+x_{\nu_1}-n}
\ \ \ts
\prod_{k=1}^{d-1}\ \ 
\frac{\Pt_{\ts i_kj_k}}{2\ts x_k-n-1}
\end{equation}
where $i$ and $i_1\com j_1\lcd i_{d-1}\com j_{d-1}$ are 
distinct indices taken from $1\lcd\nu_1-1$ 
so that different are all the sets of $d-1$ unordered pairs 
$$
\{i_1\com j_1\}\lcd\{i_{d-1}\com j_{d-1}\}\,.
$$

Replace the running index $d\ge0$ in \eqref{rthree}
by $d-1$ where $d>0\ts$.~We~get
\begin{equation}
\label{rfive}
\sum_{d>0}\,
\sum_{\substack{i_1,\ldots,i_{d-1}\\j_1,\ldots,j_{d-1}}}\, 
\sum_i\ \ \ 
\frac{\,\Pt_{\ts i\ts\nu_1}}
{x_1+x_{\nu_1}-n}
\ \ \ts
\prod_{k=1}^{d-1}\ \ 
\frac{\Pt_{\ts i_kj_k}}{2\ts x_k-n-1}
\end{equation}
with the same assumptions on the running indices as in \eqref{rfour}.
By adding up together the sums \eqref{rfour} and \eqref{rfive}, we get
\begin{equation}
\label{rsix}
\sum_{d>0}\,
\sum_{\substack{i_1,\ldots,i_{d-1}\\j_1,\ldots,j_{d-1}}}\, 
\sum_i\ \ \ 
\frac{\,\Pt_{\ts i\ts\nu_1}}
{2\ts x_d-n-1}
\ \ \ts
\prod_{k=1}^{d-1}\ \ 
\frac{\Pt_{\ts i_kj_k}}{2\ts x_k-n-1}
\end{equation}
by the equality 
$$
2\ts x_d+2\ts d=
x_1+1+x_{\nu_1}+\nu_1\ts.
$$
The sum of \eqref{rone} and \eqref{rsix} equals \eqref{aone}.
This makes the induction step.
Thus we have finished the proof of Proposition \ref{2.75}
in the case $m=1\ts$. 

Now let $m=2\ts$.
We will begin considering this case in the same way
as we did it in the proof of Proposition \ref{2.75}.
Namely, by the definition \eqref{cmu} the operator $C(\mu)$
is the ordered product of $\Rt_{\ts pq}\ts(\tts n-x_p-x_q\tts)$
over the pairs $(p\com q)$ where $1\le p<q\le\nu_1+\nu_2\ts$;
the arrangement of the pairs is reversed lexicographical.
Without changing the product,
we can rearrange these pairs as follows.
From left to right, first come the pairs $(p\com q)$ where
$\nu_1<p<q\le\nu_1+\nu_2\ts$, second come the pairs
where $1\le p\le \nu_1$ and $\nu_1<q\le\nu_1+\nu_2\ts$;
third come the pairs where $1\le p<q\le\nu_1\ts$.
Within each of the three groups,
the arrangement of the pairs $(p\com q)$ is still reversed lexicographical.

Consider the two products of $\ts\Rt_{\ts pq}\ts(\tts n-x_p-x_q\tts)\ts$,
taken over the first and the third group of pairs
$(p\com q)\ts$. If $\g_n=\so_n$ then the restriction of
each of the two products to the subspace \eqref{subtwo} is $1\ts$. 
If $\g_n=\sp_n$ then the already settled case of $m=1$ implies that 
the restrictions of the two products to \eqref{subtwo}
coincide with that of the sum
in the second displayed line in Proposition~\ref{2.75},
where $a=1$ and $a=2$ respectively.
Now consider the product of 
$\Rt_{\ts pq}\ts(\tts n-x_p-x_q\tts)$
over the second group of pairs. We have denoted this product by $Z$,
see \eqref{asreq}.
For $\g_n=\so_n\com\sp_n$ we will show that the
operator $Z$ has the same restriction to \eqref{subtwo}
as the sum in the first displayed line in Proposition \ref{2.75},
where $a=1$ and $b=2\ts$. 
This will settle the case of $m=2\ts$.

Recall 
\eqref{btwo}. 
There for any fixed $i$ we arrange 
the factors corresponding to the indices $j=1\lcd \nu_2$
from right to left. That is, we arrange
from left to right the factors corresponding to $j=\nu_2\lcd1\ts$.
The numbers $x_{\nu_1+j}$ in \eqref{btwo}
with $j=\nu_2\lcd1$ then make a sequence
increasing by $1\ts$. This is the sequence
\begin{equation}
\label{seone}
x_{\nu_1+\nu_2}\lcd x_{\nu_1+1}\,.
\end{equation}
It was only the increasing by $1$ property of \eqref{seone}
that we used to prove that the restrictions of
the operators \eqref{atwo} and \eqref{btwo} 
to the subspace \eqref{subtwo} are~equal.
Hence we can replace the sequence \eqref{seone} in this equality
by any other sequence of length $\nu_2$ that is
increasing by $1\ts$. As a replacement,
let us use the sequence 
\begin{equation*}
\label{setwo}
n-x_{\nu_1+1}\lcd n-x_{\nu_1+\nu_2}\ts.
\end{equation*}
Since $k\le\nu_1$ in \eqref{atwo}, 
then we get the equality in $(\ts\End\CC^n)^{\tts\ot\tts(\nu_1+\nu_2)}$
\begin{gather}
\label{repone}
\Bigl(\ \ 
\sum_{d\ge0}\ 
\sum_{\substack{i_1,\ldots,i_d\\j_1,\ldots,j_d}}\ \, 
\prod_{k=1}^d\ \ \ 
\frac{P_{\ts i_kj_k}}{x_k+x_{\nu_1+1}-n}
\ \,\Bigr)\,
(A_{\nu_1}\ot A_{\nu_2})\ =
\\[2pt]
\label{reptwo}
\prod_{i=1,\ldots,\nu_1}^{\longrightarrow}
\Bigl(\,\,
\prod_{j=1,\ldots,\nu_2}^{\longleftarrow}
R_{\ts i,\nu_1+j}\ts(\tts n-x_i-x_{\nu_1+\nu_2-j+1}\tts)
\Bigr)
\times
(A_{\nu_1}\ot A_{\nu_2})\,.
\end{gather}

The indices $i_1\lcd i_d$ and $j_1\lcd j_d$ in \eqref{repone}
have the same range as in the first displayed line in Proposition \ref{2.75}.
By applying to \eqref{repone} the operator 
conjugation relative to the form 
$\langle\ ,\,\rangle\ts$ in each of the first $\nu_1$ tensor factors of 
$(\ts\End\CC^n)^{\tts\ot\tts(\nu_1+\nu_2)}$ we get the product
\begin{gather}
\nonumber
(A_{\nu_1}\ot1)\,
\Bigl(\ \ 
\sum_{d\ge0}\ 
\sum_{\substack{i_1,\ldots,i_d\\j_1,\ldots,j_d}}\ \, 
\prod_{k=1}^d\ \ \ 
\frac{\Pt_{\ts i_kj_k}}{x_k+x_{\nu_1+1}-n}
\ \,\Bigr)\,
(1\ot A_{\nu_2})\,=
\\
\label{suma}
\Bigl(\ \ 
\sum_{d\ge0}\ 
\sum_{\substack{i_1,\ldots,i_d\\j_1,\ldots,j_d}}\ \, 
\prod_{k=1}^d\ \ \ 
\frac{\Pt_{\ts i_kj_k}}{x_k+x_{\nu_1+1}-n}
\ \,\Bigr)\,
(A_{\nu_1}\ot A_{\nu_2})\,.
\end{gather}
The sum over $d\ge0$ in \eqref{suma} coincides with
the sum in the first displayed line in Proposition \ref{2.75}
for $a=1$ and $b=2$, by the equalities \eqref{xk}
for $k\le\nu_1$ and 
$$
x_{\nu_1+1}=\la_2+\rho_2+(n-1)/2\,.
$$

Note that the product \eqref{repone} commutes 
with the operator on $(\CC^n)^{\tts\ot\tts(\nu_1+\nu_2)}$
reversing the order of the last $\nu_2$ tensor factors.
Let us conjugate the product \eqref{reptwo} by this operator. 
The result is the product
\begin{gather*}
\prod_{i=1,\ldots,\nu_1}^{\longrightarrow}
\Bigl(\,\,
\prod_{j=1,\ldots,\nu_2}^{\longleftarrow}
R_{\ts i,\nu_1+\nu_2-j+1}\ts(\tts n-x_i-x_{\nu_1+\nu_2-j+1}\tts)
\Bigr)
\times
(A_{\nu_1}\ot A_{\nu_2})
\\[8pt]
=\ 
\prod_{i=1,\ldots,\nu_1}^{\longrightarrow}
\Bigl(\,\,
\prod_{j=1,\ldots,\nu_2}^{\longrightarrow}
R_{\ts i,\nu_1+j}\ts(\tts n-x_i-x_{\nu_1+j}\tts)
\Bigr)
\times
(A_{\nu_1}\ot A_{\nu_2})\ =
\\[8pt]
(A_{\nu_1}\ot 1)
\times
\prod_{i=1,\ldots,\nu_1}^{\longleftarrow}
\Bigl(\,\,
\prod_{j=1,\ldots,\nu_2}^{\longrightarrow}
R_{\ts i,\nu_1+j}\ts(\tts n-x_i-x_{\nu_1+j}\tts)
\Bigr)
\times
(1\ot A_{\nu_2})\,.
\end{gather*}
By applying to 
the last displayed line the operator conjugation relative to
$\langle\ ,\,\rangle\ts$ in each of the first $\nu_1$ tensor factors of 
$(\ts\End\CC^n)^{\tts\ot\tts(\nu_1+\nu_2)}$ we get the product
\begin{equation}
\label{laslin}
Z
\times
(A_{\nu_1}\ot A_{\nu_2})\,.
\end{equation}
The equality of \eqref{repone} and \eqref{reptwo}
implies the equality of \eqref{suma} and \eqref{laslin}.
The latter equality settles the case of $m=2\ts$.
\qed
\end{proof}

The operators $B(\mu)$ and $C(\mu)$ preserve the subspace 
\eqref{subspace}.
Due to Lemmas \ref{2.1}\tts,\tts{2.2} and 
Propositions \ref{2.3}\tts,\tts\ref{2.6} the restriction 
of operator $B(\mu)\,C(\mu)$ to this subspace 
is an $\Y(\g_n)\ts$-intertwining operator \eqref{intoper2}.
Let $J(\mu)$ be this restriction
divided by the rational functions \eqref{beigen} and \eqref{ceigen},
and in the case of $\g_n=\sp_n$
also divided by the rational function \eqref{deigen}.
Then by Propositions \ref{2.4} and \ref{2.7}
\begin{equation}
\label{jmu}
J(\mu):\ts
\ph_{\nu_1}\ot\ts\ldots\ts\ot\ph_{\nu_m}\mapsto\ts
\ph_{\nu_1}\ot\ts\ldots\ts\ot\ph_{\nu_m}\ts.
\end{equation}

\begin{theorem*}
\label{2.8}
For a fixed weight $\nu=\la-\mu+\ka$
the rational function $J(\mu)$ is regular at any
point $\mu\in\h_m^*$ where the weight $\la+\rho$ is dominant.
\end{theorem*}

The operator valued rational function $J(\mu)$
does not vanish at any point $\mu\in\h_m^*$ due to the
normalization \eqref{jmu}. 
The regularity of $J(\mu)$ was proved in \cite{KN}
for all $\mu$ where the pair $(\la\com\mu)$ is good\tts;
our Theorem \ref{2.8} is more general. 

Below we will give explicit formulas for 
the operator $J(\mu)$ with dominant $\la+\rho\ts$.
Theorem \ref{2.8} will follow from these formulas.
For related results see the work of Isaev and Molev \cite{IM}. 
Note that both Proposition~\ref{2.75} and Theorem~\ref{2.8}
can also be proved by
using the arguments from \cite[Subsections 1.4 and 4.4]{KN}.


\bigskip\noindent
\textbf{3.4.}
Let the weight $\la+\rho$ of $\f_m=\so_{2m}$ or of $\f_m=\sp_{2m}$
be dominant. Then $\la_1\lcd\la_m$ obey the inequalities 
\eqref{domone} where $\rho_1\lcd\rho_m$ are now the labels of the half-sum
of positive roots of $\f_m\ts$. But for any given 
$a<b$ the difference $\rho_a-\rho_b=b-a$ is now the same as it was for
$\gl_m\ts$. Moreover, for $\f_m$ we have the same relation
$$
\la_a-\la_b=\mu_a-\mu_b+\nu_a-\nu_b
$$
as we had for $\gl_m\ts$. Any of our 
two proofs of Theorem \ref{2.5} now shows that
the operator $B(\mu)$ divided by \eqref{beigen}
has a regular restriction to the subspace~\eqref{subspace}. 
Moreover, each of the two
proofs gives an explicit formula for the restriction. 

We will give two parallel proofs of the regularity  
of restriction to \eqref{subspace}
of the operator $C(\mu)$
divided by \eqref{ceigen}, and in the case
$\f_m=\sp_{2m}$ also divided by \eqref{deigen}.
We will keep assuming that the weight $\la+\rho$ is dominant.
In particular,
\begin{align}
\label{domtwo}
\la_a+\la_b+\rho_a+\rho_b\,&\neq\,1\com2\com\,\ldots
\quad\text{whenever}\quad
a<b
\\[2pt]
\label{domthree}
\la_a+\rho_a\,&\neq\,1\com2\com\,\ldots
\quad\text{if}\quad
\f_m=\sp_{2m}\ts.
\end{align}
Then the operator $C(\mu)$
has a regular restriction to \eqref{subspace}
by Proposition \ref{2.75}.
Let us now prove the last fact directly, that is without
using Proposition \ref{2.75}.

Take any pair of indices $(p\com q)$ with $p<q$
and consider the corresponding factor $\Rt_{\ts pq}\ts(\tts n-x_p-x_q\tts)$
of the product \eqref{cmu}. Let $a$ and $b$ be the numbers of the segments
of the sequence $1\lcd N$ containing the indices $p$ and $q$ respectively. 
Let $i$ and $j$ be the numbers of the indices $p$ and $q$
within their segments. Then
\begin{equation}
\label{xpqn}
x_p+x_q-n=
\la_a+\rho_a+\la_b+\rho_b-i-j+1\,.
\end{equation}

First, suppose that $p$ and $q$ belong to different segments,
so that $a<b\ts$. Then the right hand side of the equality \eqref{xpqn} 
is not zero by \eqref{domtwo}.
Therefore the factor $\Rt_{\ts pq}\ts(\tts n-x_p-x_q\tts)$ of 
the product \eqref{cmu} with $a<b$ is regular.

Next, suppose that $p$ and $q$ belong to the same segment, that is $a=b\ts$.
If $\f_m=\so_{2m}$ then the factor $\Rt_{\ts pq}\ts(\tts n-x_p-x_q\tts)$ can
be skipped without changing the restriction of \eqref{cmu} to the subspace
\eqref{subspace}. If $\f_m=\sp_{2m}$ then by
\eqref{domthree} the right hand side of the equality
\eqref{xpqn} is not zero whenever the sum 
$i+j$ is odd. In particular, if $\f_m=\sp_{2m}$ then
the factor $\Rt_{\ts pq}\ts(\tts n-x_p-x_q\tts)$ 
with $q=p+1$ is regular, because for this factor $j=i+1\ts$.

Last, suppose that $q>p+1$ while $a=b\ts$.
Then the pair following
$(p\com q)$ in the reversed lexicographical ordering is 
$(p\com q-1)\ts$. Moreover, then the index $q-1$ belongs
to the same segment as $p$ and $q\ts$.
Take the product of the
factors in \eqref{cmu} corresponding the pairs
following $(p\com q-1)\ts$. Multiply this product by
\begin{equation}
\label{pq}
R_{q-1,q}\ts(x_{q-1}-x_q)=R_{q-1,q}\ts(1)=1-P_{q-1,q}
\end{equation}
on the right.
By using the relation \eqref{cybe} repeatedly,
one shows that the resulting product
is also divisible by \eqref{pq} on the left.
Due to \eqref{nopolet} 
we can then replace in \eqref{cmu} 
the product of two adjacent factors 
$$
\Rt_{pq}\ts(n-x_p-x_q)\,
\Rt_{p,q-1}\ts(n-x_p-x_{q-1})\,
$$
by
$$
1-(\ts\Pt_{pq}+\Pt_{p,q-1}\ts)/(n-x_p-x_{q-1})
$$
without changing the restriction of
the operator \eqref{cmu} to the subspace   
\eqref{subspace}. But the replacement
is regular at $x_p+x_q=n\ts$.
This observation completes our second proof of the regularity
of the restriction of $C(\mu)$ to the subspace \eqref{subspace},
whenever the weight $\la+\rho$ of $\f_m$ is dominant.

Now recall that defining the operator $J(\mu)$
involves dividing $C(\mu)$ by \eqref{ceigen}, and also
dividing by \eqref{deigen} in the case $\f_m=\sp_{2m}\ts$.
So we have to consider the zeroes of the rational functions
\eqref{ceigen} and \eqref{deigen}. In Subsection 3.5
for $m=1$ and $\g_n=\sp_n$ we will prove that 
\eqref{cone} annihilates the subspace \eqref{lanone} whenever
\begin{equation}
\label{condone}
\la_1+\rho_1=0
\quad\text{and}\quad
2\ts\nu_1>n\ts.
\end{equation}
In Subsection 3.6 for $m=2$ and both $\g_n=\so_n\com\sp_n$ 
we will prove that the operator
\eqref{asreq} annihilates the subspace \eqref{subtwo} whenever
\begin{equation}
\label{condtwo}
\la_1+\la_2+\rho_1+\rho_2=0
\quad\text{and}\quad
\nu_1+\nu_2>n\ts.
\end{equation}
Theorem \ref{2.8} for any $m\ge1$ will then follow
from the definitions \eqref{ceigen}~and~\eqref{deigen}. 


\bigskip\noindent
\textbf{3.5.}
For $m=1$ and $\g_n=\sp_n$
consider the operator \eqref{cone} on 
the vector space $(\CC^n)^{\tts\ot\tts\nu_1}\ts$.
Here the positive integer $n$ is even. 
For $p=1\lcd\nu_1$ introduce the rational function of $x\in\CC$
taking values in the operator algebra $(\ts\End\CC^n)^{\tts\ot\tts\nu_1}$
\begin{equation}
\label{dxl}
D(x\com p\tts)\,=\!
\prod_{j=2,\ldots,p}^{\longleftarrow}
\!\Bigl(\ 
\prod_{i=1,\ldots,j-1}^{\longleftarrow}
\Rt_{ij}\ts(\tts i+j-2\tts x-1\tts)
\Bigr)\,.
\end{equation}
By \eqref{xk} we have
\begin{equation}
\label{cd}
C(\mu)=D(\tts\la_1+\rho_1\com\nu_1)
\end{equation}
while $D(x\com 1\tts)=1$ by definition. Put $D(x\com 0\tts)=1\ts$.
Like we did at the beginning of the proof of Proposition \ref{2.75},
let us relate two operators on the vector space 
$(\CC^n)^{\tts\ot\tts\nu_1}$
by the symbol~$\,\equiv$ if their actions on
the subspace \eqref{lanone} coincide.

\begin{lemma*}
\label{3.55}
For each $p=2\lcd\nu_1$ we have the relation in 
$(\ts\End\CC^n)^{\tts\ot\tts\nu_1}$
\begin{equation}
\label{35}
\Pt_{\,p-1,p}\,D(x\com p\tts)\,\equiv\,
\frac{x+n/2-p+1}{x-1}\ 
\Pt_{\,p-1,p}\,D(x-1\com p-2\tts)\,.
\end{equation}
\end{lemma*}

\begin{proof}
Using Proposition \ref{2.75} where $m=1$ whereas
the numbers
$\la_1+\rho_1$ and $\nu_1$ are replaced by $x$ and $p$ respectively,
we obtain the relation
\begin{equation}
\label{dxnu}
D(x\com p)\,\equiv\,
\displaystyle
\sum_{d\ge 0}\ 
\sum_{\substack{i_1,\ldots,i_d\\j_1,\ldots,j_d}}\ \, 
\prod_{k=1}^d\ \ \, 
\frac{\Pt_{\ts i_kj_k}}{2\ts(\ts x-k\ts)}
\end{equation}
where $\,i_1\com j_1\lcd i_d\com j_d\ts$ are pairwise 
distinct indices from the sequence $1\lcd p$ taken so that
different are all the sets of $d$ unordered pairs \eqref{dup}.
Like in the proof of Proposition \ref{2.75}
we assume that $1$ is the only term in \eqref{dxnu} with $d=0\ts$.

For any $d\ge0$ and for any 
choice of the set of $d$ pairs \eqref{dup} made as above
consider the corresponding product 
$\Pt_{\ts i_1j_1}\ldots\Pt_{\ts i_dj_d}$
showing in \eqref{dxnu}. Multiply
the product by $\Pt_{\ts p-1,p}$ on the left.
If neither of the indices $p-1\com p$ occurs 
in the pairs \eqref{dup} then leave the result
of multiplication as it is,
$
\,\Pt_{\ts p-1,p}\,\Pt_{\ts i_1j_1}\ldots\Pt_{\ts i_dj_d}\,.
$

Next, suppose that exactly one of the indices $p-1\com p$ 
occurs in \eqref{dup}.
We can assume that then $j_d=p-1$ or $j_d=p$ without further 
loss of generality. Put $j=p$ or $j=p-1$ respectively. Then
\begin{align}
\nonumber
\Pt_{\ts p-1,p}\,\Pt_{\ts i_1j_1}\ldots\Pt_{\ts i_dj_d}
&=\Pt_{\ts p-1,p}\,\Pt_{\ts i_1j_1}\ldots \Pt_{\ts i_{d-1}j_{d-1}}
\,P_{\ts i_d\ts j}\,\equiv
\\
\label{termone}
&-\Pt_{\ts p-1,p}\,\Pt_{\ts i_1j_1}\ldots \Pt_{\ts i_{d-1}j_{d-1}}\,.
\end{align}
Note that in either case, that is $j_d=p-1$ or $j_d=p\ts$, 
for any given distinct indices 
$i_1\com j_1\lcd i_{d-1}\com j_{d-1}$
taken from $1\lcd p-2$ there are exactly $p-2\tts d$
choices of the index $i_d$ yielding the same term \eqref{termone}
where $i_d$ does not occur. Counting both cases, we will get 
the term \eqref{termone} with multiplicity $2\ts(\tts p-2\tts d\ts)\ts$.

Finally, suppose that both of the indices $p-1\com p$ occur 
in \eqref{dup}. If they occur in the same pair, then without
further loss of generality we may assume that
$i_d=p-1$ and $j_d=p\ts$. Then
\begin{equation}
\label{termn}
\Pt_{\ts p-1,p}\,\Pt_{\ts i_1j_1}\ldots\Pt_{\ts i_dj_d}
\,=\,n\,\Pt_{\ts p-1,p}\,\Pt_{\ts i_1j_1}\ldots \Pt_{\ts i_{d-1}j_{d-1}}\,.
\end{equation}
If $p-1\com p$ occur in different pairs in \eqref{dup} then without
further loss of generality we may assume that
$j_{d-1}=p-1$ and $i_d=p\ts$. Then
\begin{align}
\nonumber
\Pt_{\ts p-1,p}\,\Pt_{\ts i_1j_1}\ldots\Pt_{\ts i_dj_d}
&=\Pt_{\ts p-1,p}\,\Pt_{\ts i_1j_1}\ldots\Pt_{\ts i_{d-2}j_{d-2}}\,
P_{\ts i_{d-1},p}\,\Pt_{\ts p,j_d}\,\equiv
\\
\label{termtwo}
&-\Pt_{\ts p-1,p}\,\Pt_{\ts i_1j_1}\ldots\Pt_{\ts i_{d-2}j_{d-2}}\,
\Pt_{\ts i_{d-1}j_d}\,.
\end{align}
Without altering the value of the term \eqref{termtwo}, 
we can either exchange the pair $(i_{d-1}\com j_d)$ with 
any of the pairs $(i_1\com j_1)\lcd (i_{d-2}\com j_{d-2})$
or swap the two indices $i_{d-1}\com j_d$ between each other.
Counting these together with the initial choice of $i_{d-1}\com j_d$ 
we will get the term \eqref{termtwo} with multiplicity 
$2\ts(\tts d-1\tts)\ts$. Note that here the term \eqref{termtwo}
arises only when $d\ge2$. 
But for $d=1$ we have $2\ts(\tts d-1\tts)=0\ts$.

Let us multiply the relation \eqref{dxnu}
by $\Pt_{\ts p-1,p}$ on the left, and then perform the
summation over the indices $d$ and $\,i_1\com j_1\lcd i_d\com j_d\ts$.
The terms \eqref{termone},\eqref{termn} occur only when $d\ge1\ts$.
Now replace the index $d-1$ in \eqref{termone},\eqref{termn} by $d\ge 0\ts$. 
After this replacement, the multiplicity of \eqref{termone} will become
$\,2\ts(\tts p-2\tts d-2\ts)\ts$.

Further, rename the running index $j_d$ in \eqref{termtwo}
by $j_{d-1}$ and then replace $d-1$ by $d\ts$.
After this, the corresponding multiplicity is $2\ts d\ts$.
Then $d\ge1\ts$. But we can also sum over $d\ge0\ts$, because
the multiplicity in the case $d=0$ is zero.

After these replacements, the product
$\Pt_{\ts p-1,p}\,D(x\com p\tts)$ gets related by $\equiv$ to
\begin{gather*}
\sum_{d\ge 0}\ 
\sum_{\substack{i_1,\ldots,i_d\\j_1,\ldots,j_d}}
\Pt_{\,p-1,p}\,
\biggl(\,
1+\frac{n-2\ts(\ts p-2\tts d-2\ts)-2\tts d}{2\ts(\ts x-d-1\ts)}
\,\biggr)\ 
\prod_{k=1}^d\  
\frac{\Pt_{\ts i_kj_k}}{2\ts(\ts x-k\ts)}
\ =
\\
\sum_{d\ge 0}\ 
\sum_{\substack{i_1,\ldots,i_d\\j_1,\ldots,j_d}}
\Pt_{\ts p-1,p}\,
\ \frac{x+n/2-p+1}{x-d-1}\ 
\prod_{k=1}^d\  
\frac{\Pt_{\ts i_kj_k}}{2\ts(\ts x-k\ts)}
\ =
\\
\sum_{d\ge 0}\ 
\sum_{\substack{i_1,\ldots,i_d\\j_1,\ldots,j_d}}
\Pt_{\ts p-1,p}\,
\ \frac{x+n/2-p+1}{x-1}\ 
\prod_{k=1}^d\  
\frac{\Pt_{\ts i_kj_k}}{2\ts(\ts x-k-1\ts)}
\end{gather*}
where $\,i_1\com j_1\lcd i_d\com j_d\ts$ are distinct
indices taken from the sequence $1\lcd p-2$ such that
different are all the sets of $d$ unordered pairs \eqref{dup}.
But the sum in the last displayed line is related by $\,\equiv\,$
to 
the right hand side of \eqref{35}. Here we use the relation
\eqref{dxnu} with $x-1$ and $p-2$ instead of $x$ and $p$ respectively. 
\qed
\end{proof}

From now until the end of this subsection
we assume that $2\ts\nu_1>n\ts$. Denote 
$$
l=\nu_1-n/2\ts.
$$
It is a classical fact \cite[Section VI.3]{W}
that then the subspace \eqref{lanone}
is contained in the span of the images 
of all operators on $(\CC^n)^{\tts\ot\tts\nu_1}$ of the form
$
\Pt_{\ts i_1j_1}\ldots\Pt_{\ts i_lj_l}
$
where $i_1\com j_1\lcd i_l\com j_l$
are any pairwise distinct indices
taken from $1\lcd\nu_1\ts$.
To show that $C(\mu)\equiv0$ under the conditions 
\eqref{condone}, it now suffices to prove
\begin{equation}
\label{pca}
\Pt_{\ts i_1j_1}\ldots\Pt_{\ts i_lj_l}\,C(\mu)\,\equiv\,0 
\quad\text{for}\quad
\la_1+\rho_1=0
\end{equation}
and for all those $i_1\com j_1\lcd i_l\com j_l\ts$.
Here we also use the equality $\Pt^{\,2}=n\,\Pt$
and the fact that the operator 
$C(\mu)$ with generic $\mu$ preserves the subspace \eqref{lanone}.
Furthermore, due to the latter fact,
it suffices to prove the relation \eqref{pca} for any single
choice of the indices $i_1\com j_1\lcd i_l\com j_l\ts$. Let us choose
$$
i_1=\nu_1-2\ts l+1\com\ j_1=
\nu_1-2\ts l+2\com\ \ldots\ ,\ 
i_l=\nu_1-1\com\ j_l=\nu_1\ts.
$$
By using \eqref{cd} and then applying Lemma \ref{3.55} repeatedly,
namely applying it $l$ times, we get the following relation between 
rational functions of $\mu\in\h_1^*\,$:  
\begin{align*}
&\Pt_{\ts\nu_1-2l+1\tts,\ts\nu_1-2l+2}
\,\ldots\ts
\Pt_{\ts\nu_1-1,\tts\nu_1}\,C(\mu)
\,\equiv\,
\frac{\,\la_1+\rho_1}{\mu_1+\rho_1}\ \times
\\
&\Pt_{\ts\nu_1-2l+1\tts,\ts\nu_1-2l+2}
\,\ldots\ts
\Pt_{\ts\nu_1-1,\tts\nu_1}\,
D(\ts\mu_1+\rho_1\com n-\nu_1\ts)\,. 
\end{align*}

To get the latter relation, we also used the equality
$\nu_1-2\ts l=n-\nu_1\ts$.
If $\la_1+\rho_1=0$ then $\mu_1+\rho_1=-\,l\ts$, and
the fraction in the above display equals zero.
The last factor in that display
is regular at $\la_1+\rho_1=0$ by the definition \eqref{dxl}.
This proves \eqref{pca} for our 
choice of $i_1\com j_1\lcd i_l\com j_l\ts$. 

As a rational function of $\mu\ts$, the last displayed
factor can be replaced by the sum at the right hand side
of the relation~\eqref{dxnu}
where $x=\mu_1+\rho_1$ and $p=n-\nu_1\ts$.
Arguing like in Subsection 3.4, we can also provide
a multiplicative formula
for the value of that factor whenever 
$\mu_1+\rho_1\neq1\com2\com\ts\ldots\,\,$. 

Multiplying the last displayed relation on the left
by operators on $(\CC^n)^{\tts\ot\tts\nu_1}$ which permute the
$\nu_1$ tensor factors, 
we obtain analogues of that relation
for all other choices of 
$i_1\com j_1\lcd i_l\com j_l\,$.
Therefore in the case $m=1$ and $2\ts\nu_1>n\ts$, 
the last displayed relation determines the action on
the subspace \eqref{lanone} for~any value
of the function $C(\mu)$ divided by~\eqref{deigen},
whenever $\mu_1+\rho_1\neq1\com2\com\ts\ldots\,\,$. 


\bigskip\noindent
\textbf{3.6.}
For $m=2$ and $\g_n=\so_n\com\sp_n$
consider the operator \eqref{asreq} on
$(\CC^n)^{\tts\ot\tts(\nu_1+\nu_2)}\ts$.
For $p=1\lcd\nu_1$ and $q=\nu_1+1\lcd\nu_1+\nu_2$
introduce the rational function of $x\in\CC$
taking values in the algebra $(\ts\End\CC^n)^{\tts\ot\tts(\nu_1+\nu_2)}$
\begin{equation}
\label{dy}
D(x\com p\com q\tts)\,=\!
\prod_{i=1,\ldots,p}^{\longleftarrow}
\!\Bigl(\ 
\prod_{j=1,\ldots,q-\nu_1}^{\longleftarrow}
\Rt_{ij}\ts(\tts i+j-x-1\tts)
\Bigr)\,.
\end{equation}
By \eqref{xk}
\begin{equation}
\label{zd}
Z=D(\tts\la_1+\la_2+\rho_1+\rho_2\com\ts\nu_1\com\nu_1+\nu_2)\,.
\end{equation}
Put $D(x\com p\com \nu_1\tts)=D(x\com 0\com q\tts)=1\ts$.
Let us relate operators on $(\CC^n)^{\tts\ot\tts(\nu_1+\nu_2)}$
by the symbol~$\,\equiv$ if their actions coincide on
the subspace \eqref{subtwo}.

\begin{lemma*}
\label{3.65}
For any $\,p=1\lcd\nu_1\,$ and $\,q=\nu_1+1\lcd\nu_1+\nu_2\,$
we have
\begin{equation}
\label{36}
\Pt_{pq}\,D(x\com p\com q\tts)\,\equiv\,
\frac{x+n-p-q+\nu_1+1}{x-1}\ 
\Pt_{pq}\,D(\tts x-1\com p-1\com q-1\tts)\,.
\end{equation}
\end{lemma*}

\begin{proof}
At the end of the proof of Proposition \ref{2.75}
we established the equality of the expressions
\eqref{suma} and \eqref{laslin}. 
Replacing the numbers $\la_1+\la_2+\rho_1+\rho_2$ and $\nu_1\com\nu_2$
in that equality by $x$ and $p\com q-\nu_1$ respectively,
we get the relation
\begin{equation}
\label{dynu}
D(x\com p\com q)\,\equiv\,
\displaystyle
\sum_{d\ge 0}\ 
\sum_{\substack{i_1,\ldots,i_d\\j_1,\ldots,j_d}}\ \, 
\prod_{k=1}^d\ \ \, 
\frac{\,\Pt_{\ts i_kj_k}}{x-k}
\end{equation}
where $i_1\lcd i_d$ and $j_1\lcd j_d$ are
distinct indices taken respectively from the sequences
$1\lcd p$ and $\nu_1+1\lcd q$ so that 
different are the corresponding sets of $d$ pairs \eqref{dp}. 
We assume that $1$ is the only summand in \eqref{dynu} with $d=0\ts$.

For any $d\ge0$ and for any 
choice of the set of $d$ pairs \eqref{dp} made as above
consider the corresponding product 
$\Pt_{\ts i_1j_1}\ldots\Pt_{\ts i_dj_d}$
showing in \eqref{dynu}. Multiply
the product by $\Pt_{pq}$ on the left.
If neither of the indices $p\com q$ occurs 
in the pairs \eqref{dp} then leave the result
of multiplication as it is, that is
$
\,\Pt_{pq}\,\Pt_{\ts i_1j_1}\ldots\Pt_{\ts i_dj_d}\,.
$

Next, suppose that exactly one of the indices $p\com q$ 
occurs in \eqref{dp}.
We can assume that then $i_d=p$ or $j_d=q$ without further 
loss of generality. Then
\begin{equation}
\label{termonepq}
\Pt_{pq}\,\Pt_{\ts i_1j_1}\ldots\Pt_{\ts i_dj_d}
\equiv
-\,\Pt_{pq}\,\Pt_{\ts i_1j_1}\ldots \Pt_{\ts i_{d-1}j_{d-1}}\,.
\end{equation}
In the first case, that is if $i_d=p\ts$,
for any given $i_1\com j_1\lcd i_{d-1}\com j_{d-1}$
there are exactly $q-\nu_1-d$
choices of the index $j_d$ yielding the 
same right hand side of \eqref{termonepq},
where $j_d$ does not occur. 
In the second case, that is if $j_d=q\ts$,
for any given $i_1\com j_1\lcd i_{d-1}\com j_{d-1}$
there are exactly $p-d$
choices of the index $i_d$ yielding the same 
right hand side of \eqref{termonepq},
where $i_d$ does not occur. 
Counting both cases, we will get 
the term \eqref{termonepq} with multiplicity 
$p+q-\nu_1-2\tts d\ts$.

Last suppose $p\com q$ both occur 
in \eqref{dp}. If they occur in the same pair, then without
further loss of generality we may assume that
$i_d=p$ and $j_d=q\ts$.~Then
\begin{equation}
\label{termnpq}
\Pt_{pq}\,\Pt_{\ts i_1j_1}\ldots\Pt_{\ts i_dj_d}
\,=\,n\,\Pt_{\ts pq}\,\Pt_{\ts i_1j_1}\ldots \Pt_{\ts i_{d-1}j_{d-1}}\,.
\end{equation}
If $p\com q$ occur in different pairs in \eqref{dp} then without
further loss of generality we may assume that
$j_{d-1}=q$ and $i_d=p\ts$. Then
\begin{equation}
\label{termtwopq}
\Pt_{pq}\,\Pt_{\ts i_1j_1}\ldots\Pt_{\ts i_dj_d}\equiv
-\,\Pt_{pq}\,\Pt_{\ts i_1j_1}\ldots\Pt_{\ts i_{d-2}j_{d-2}}\,
\Pt_{\ts i_{d-1}j_d}\,.
\end{equation}
Without altering the product at the right hand side of \eqref{termtwopq}, 
we can exchange the pair $(i_{d-1}\com j_d)$ with 
any of the pairs $(i_1\com j_1)\lcd (i_{d-2}\com j_{d-2})\ts$.
Counting these together with the initial choice of $(i_{d-1}\com j_d)$ 
we will get the term \eqref{termtwo} with multiplicity 
$d-1\ts$. Note that here the term \eqref{termtwopq} arises only when $d\ge2$. 

Let us multiply the relation \eqref{dynu}
by $\Pt_{pq}$ on the left, and then perform the
summation over the indices $d$ and $\,i_1\com j_1\lcd i_d\com j_d\ts$.
The terms \eqref{termonepq},\eqref{termnpq} occur only when $d\ge1\ts$. Now
replace the index $d-1$ in \eqref{termonepq},\eqref{termnpq} by $d\ge 0\ts$. 
After this replacement, the multiplicity of \eqref{termonepq} becomes
$p+q-\nu_1-2\tts d-2\ts$.

Further, let us rename the running index $j_d$ at the right hand 
side of \eqref{termtwopq}
by $j_{d-1}\ts$, and then replace $d-1$ by $d$ there.
Having done this, the corresponding multiplicity becomes $d\ts$.
Now $d\ge1\ts$. But we can also sum over $d\ge0\ts$, because
the multiplicity in the case $d=0$ is zero.

After these replacements, the product
$\Pt_{pq}\,D(x\com p\com q\tts)$ gets related by $\equiv$ to
\begin{gather*}
\sum_{d\ge 0}\ 
\sum_{\substack{i_1,\ldots,i_d\\j_1,\ldots,j_d}}
\Pt_{pq}\,
\biggl(\,
1+\frac{n-(\ts p+q-\nu_1-2\tts d-2\ts)-d}{x-d-1}
\,\biggr)\ 
\prod_{k=1}^d\  
\frac{\,\Pt_{\ts i_kj_k}}{x-k}
\ =
\\
\sum_{d\ge 0}\ 
\sum_{\substack{i_1,\ldots,i_d\\j_1,\ldots,j_d}}
\Pt_{pq}\,
\ \frac{x+n-p-q+\nu_1+1}{x-1}\ 
\prod_{k=1}^d\  
\frac{\Pt_{\ts i_kj_k}}{x-k-1}
\end{gather*}
where $i_1\lcd i_d$ and $j_1\lcd j_d$ are pairwise
distinct indices taken respectively from 
$1\lcd p-1$ and $\nu_1+1\lcd q-1$ so that 
different are the corresponding sets of $d$ pairs \eqref{dp}. 
The sum in the last displayed line is related by $\,\equiv\,$
to 
the right hand side of \eqref{36}. Here we use the relation \eqref{dynu}
with $x-1\com p-1\com q-1$ instead of $x\com p\com q$ respectively. 
\qed
\end{proof}

From now until the end of this subsection
assume that $\nu_1+\nu_2>n\ts$. Denote
$$
l=\nu_1+\nu_2-n\ts.
$$
Then the subspace \eqref{subtwo}
is contained in the span of the images 
of all operators on $(\CC^n)^{\tts\ot\tts(\nu_1+\nu_2)}$ of the form
$
\Pt_{\ts i_1j_1}\ldots\Pt_{\ts i_lj_l}
$
where $i_1\lcd i_l$ and $j_1\lcd j_l$ are
pairwise distinct indices 
from the sequences $1\lcd\nu_1$ and $\nu_1+1\lcd\nu_1+\nu_2$ 
respectively\tts;
see for instance the proof of \cite[Lemma V.7.B]{W}.
To show that $Z\equiv0$ under the conditions \eqref{condtwo},
it therefore suffices to prove that
\begin{equation}
\label{pza}
\Pt_{\ts i_1j_1}\ldots\Pt_{\ts i_lj_l}\,Z\,\equiv\,0 
\quad\text{for}\quad
\la_1+\la_2+\rho_1+\rho_2=0
\end{equation}
and for all those $i_1\com j_1\lcd i_l\com j_l\ts$.
Here we also use the 
fact that for any generic $\mu\in\h_2^*$ the operator 
$Z$  preserves the subspace \eqref{subtwo}.
Due to the latter fact,
it suffices to prove \eqref{pza} for any single choice of 
$i_1\com j_1\lcd i_l\com j_l\ts$. Let us choose
$$
i_1=\nu_1-l+1\com\ 
j_1=\nu_1+\nu_2-l+1\com\ \ldots\ ,\ 
i_l=\nu_1\com\ j_l=\nu_1+\nu_2\ts.
$$
By using \eqref{zd} and then applying Lemma \ref{3.65} repeatedly,
namely applying it $l$ times, we get the following relation between 
rational functions of $\mu\in\h_2^*\,$:  
\begin{align*}
&\Pt_{\ts\nu_1-l+1\tts,\nu_1+\nu_2-l+1}
\,\ldots\ts
\Pt_{\ts\nu_1,\tts\nu_1+\nu_2}\,Z
\,\equiv\,
\frac{\,\la_1+\la_2+\rho_1+\rho_2}{\mu_1+\mu_2+\rho_1+\rho_2}\ \times
\\
&\Pt_{\ts\nu_1-l+1\tts,\nu_1+\nu_2-l+1}
\,\ldots\ts
\Pt_{\ts\nu_1,\tts\nu_1+\nu_2}\,
D(\ts\mu_1+\mu_2+\rho_1+\rho_2\com n-\nu_2\com n\ts)\,. 
\end{align*}

If $\la_1+\la_2+\rho_1+\rho_2=0$ then $\mu_1+\mu_2+\rho_1+\rho_2=-\,l\ts$,
and the fraction in the above display equals zero.
But the last factor in that display
is regular at $\la_1+\la_2+\rho_1+\rho_2=0$ by the definition \eqref{dy}.
This proves the relation \eqref{pza} for our particular 
choice of the indices $i_1\com j_1\lcd i_l\com j_l\ts$. 

As rational function of $\mu\ts$, the last 
factor can be replaced by the sum at right hand side
of the relation~\eqref{dynu}
where $x=\mu_1+\mu_2+\rho_1+\rho_2$ while $p=n-\nu_2$ and $q=n\ts$.
Arguing like in Subsection 3.4, we can also provide
a multiplicative formula
for the value of that factor whenever 
$\mu_1+\mu_2+\rho_1+\rho_2\neq1\com2\com\ts\ldots\,\,$. 

Further, multiplying the last displayed relation on the left
by operators on $(\CC^n)^{\tts\ot\tts(\nu_1+\nu_2)}$ which permute the
first $\nu_1$ tensor factors between themselves, 
and also permute the last $\nu_2$ tensor factors, 
we get analogues of that relation
for all other choices of 
$i_1\com j_1\lcd i_l\com j_l\,$.
So in the case $m=2$ and $\nu_1+\nu_2>n\ts$, 
the last displayed relation determines the action on
the subspace \eqref{subtwo} for any value
of the function $Z$ divided by~\eqref{ceigen},
when $\mu_1+\mu_2+\rho_1+\rho_2\neq1\com2\com\ts\ldots\,\,$. 


\bigskip\noindent
\textbf{3.7.}
In this subsection we will generalize Theorem \ref{1.2}.
Theorem \ref{2.8} allows us to determine the intertwining operator $J(\mu)$ 
of the $\Y(\g_n)\ts$-modules \eqref{intoper2}
for any $\mu\in\h_m^*\ts$, provided $\la+\rho$ is dominant.
Our generalization of Theorem \ref{1.2} is based on the next
two lemmas. For each index $c=1\lcd m-1$ regard $s_c$ as an element
of the group $H_m\ts$. The action of $s_c$ on $\h_m$ exchanges
the basis vectors $F_{cc}$ and $F_{c+1,c+1}\ts$, leaving all  
other basis vectors of $\h_m$ fixed. 
The first of the two lemmas is an analogue of Lemma \ref{2.55} for the twisted
Yangian $\Y(\g_n)\ts$.

\begin{lemma*}
\label{3.51}
Fix $c>0$. Suppose that both $\la+\rho$ and\/ $s_c(\la+\rho)$
are dominant. Then the images of 
the intertwining operators $J(\mu)$ and $J(s_c\circ\mu)$ 
corresponding to the pairs $(\la\com\mu)$ and 
$(s_c\circ\la\com s_c\circ\mu)$ 
are equivalent as\/ $\Y(\g_n)\ts$-modules.
\end{lemma*}

\begin{proof}
Let $\muc$ and $\rhoc$ be the weights of $\gl_2$ with 
the labels $\mu_c\com\mu_{c+1}$ and
$\rho_c\com\rho_{c+1}$ respectively. Let $\lac$ be 
the weight of $\gl_2$ with labels
$\mu_c+\nu_c\com\mu_{c+1}+\nu_{c+1}\ts$.
Here
$$
\mu_c+\nu_c=\la_c+n/2
\quad\text{and}\quad
\mu_{c+1}+\nu_{c+1}=\la_{c+1}+n/2\ts.
$$
The dominance of the weights $\la+\rho$ and $s_c(\la+\rho)$ of $\f_m$ 
implies that the weights $\lac+\rhoc$ and $s_1(\lac+\rhoc)$ of $\gl_2$
are also dominant. By Theorem \ref{2.5} with $m=2\ts$, 
we get an intertwining operator of $\Y(\gl_n)\ts$-modules
\eqref{fint}. It is invertible,
see the proof of Lemma \ref{2.55}. 
Like in that proof, denote by $I$ the operator which acts 
on the tensor product of $c\,$th and
$(c+1)\,$th factors of \eqref{basprod1}
as this intertwining operator \eqref{fint}, and
which acts trivially on other $m-2$ tensor factors of \eqref{basprod1}.

Similarly, by using Theorem \ref{2.5} with $m=2$ once again,
we get an invertible intertwining operator of $\Y(\gl_n)\ts$-modules
\begin{equation}
\label{mint}
\Ph^{\,-\nu_c}_{\mu_c+\rho_c+\frac12}
\ot
\Ph^{\,-\nu_{c+1}}_{\mu_{c+1}+\rho_{c+1}+\frac12}
\to\,
\Ph^{\,-\nu_{c+1}}_{\mu_{c+1}+\rho_{c+1}+\frac12}
\ot
\Ph^{\,-\nu_c}_{\mu_c+\rho_c+\frac12}\,.
\end{equation}
Now denote by $J$ the operator which acts as \eqref{mint} 
on the tensor product of $c\,$th and
$(c+1)\,$th factors of the target $\Y(\g_n)\ts$-module in
\eqref{intoper2}, and
which acts trivially on other $m-2$ tensor factors of 
the latter module. Arguing like in the end of the proof of
Lemma \ref{2.55}, that is either performing a direct calculation,
or using the irreducibility of the source and target
$\Y(\g_n)\ts$-modules in \eqref{intoper2} for any generic weight
$\mu$ of $\f_m\ts$, we obtain the relation
$$
J\,J(\mu)=J(s_c\circ\mu)\,I 
$$
whenever $\la+\rho$ and $s_c(\la+\mu)$ are dominant.
It proves Lemma \ref{3.51}, since $I$ and $J$
are invertible and intertwine
$\Y(\g_n)\ts$-modules by restriction from $\Y(\gl_n\ts)$. 
\qed
\end{proof}

Let $s_0\in H_m$ be the element which acts on $\h_m$ by mapping
$F_{11}$ to $-\ts F_{11}\ts$, and leaves all other basis vectors fixed.
Note that in the case of $\f_m=\so_{2m}$ the element $s_0$
belongs to the extended Weyl group, not to the Weyl group proper.
Further, in this case the dominance of $s_0(\la+\rho)$
is equivalent to that of $\la+\rho\ts$.

\begin{lemma*}
\label{3.52}
Suppose that both the weights $\la+\rho$ and\/ $s_0(\la+\rho)$
are dominant. Then the images of 
the intertwining operators $J(\mu)$ and $J(s_0\circ\mu)$ 
corresponding to the pairs $(\la\com\mu)$ and 
$(s_0\circ\la\com s_0\circ\mu)$ 
are similar as\/ $\Y(\g_n)\ts$-modules.
\end{lemma*}

\begin{proof}
Let $\lac\com\muc\com\rhoc$ be the weights of $\f_1$ with
labels $\la_1\com\mu_1\com\rho_1$ respectively. The weights
$\lac+\rhoc$ and $s_0(\lac+\rhoc)$ of $\f_1$ are dominant. 
For $\f_1=\sp_2$ this means that $\la_1\notin\ZZ\,\backslash\ts\{1\}\ts$.
For $\f_1=\so_2$ any weight is dominant.
By using Theorem~\ref{2.8} with $m=1$  
we get the intertwining operators of $\Y(\g_n)\ts$-modules
\begin{align}
\label{hint}
\Ph^{\,\nu_1}_{\mu_1+\rho_1+\frac12}
&\to\,
\Ph^{\,-\nu_1}_{\mu_1+\rho_1+\frac12}\,,
\\[4pt]
\label{bint}
\Ph^{\,n-\nu_1}_{\frac12-\mu_1-\rho_1}
&\to\,
\Ph^{\,\nu_1-n}_{\frac12-\mu_1-\rho_1}\,.
\end{align}
For $\f_1=\so_{2}$ each of these two operators acts as the identity,
see remarks at the beginning of Subsection 3.2.
For $\f_1=\sp_{2}$ none of these two operators acts as the identity
in general, but they are still invertible. The latter
assertion can be proved either by direct calculation,
or by using the irreducibility of all four
$\Y(\g_n)\ts$-modules in \eqref{hint} and \eqref{bint}
for generic~$\muc\ts$, that is for $\mu_1\notin\ZZ/2$.

For instance, let us prove the invertibility of 
\eqref{hint}. By applying Lemma \ref{1.15} to both the source and target
$\Y(\gl_n)\ts$-modules in \eqref{bint} 
we get an intertwining operator 
\begin{equation}
\label{dint}
\Ph^{\,-\nu_1}_{\mu_1+\rho_1+\frac12}
\to\,
\Ph^{\,\nu_1}_{\mu_1+\rho_1+\frac12}
\end{equation}
of $\Y(\sp_n)\ts$-modules.
This operator maps $\ph_{\nu_1}\mapsto\ph_{\nu_1}$
as does the operator \eqref{hint}.
Hence the operators \eqref{hint} and \eqref{dint}
are inverse to each other for $\mu_1\notin\ZZ/2$,
and therefore for any $\mu_1\in\CC$ such that  
$\la_1\notin\ZZ\,\backslash\ts\{1\}\ts$.

By using Lemma \ref{1.15} with $k=\nu_1$ and $t=\mu_1+\rho_1+\frac12$
we get an invertible intertwining operator
\begin{equation}
\label{lint}
\Ph^{\,-\nu_1}_{\mu_1+\rho_1+\frac12}
\to\,
\Phd^{\,n-\nu_1}_{\frac12-\mu_1-\rho_1}
\end{equation}
of $\Y(\gl_n)\ts$-modules. It is also
an intertwiner of $\Y(\g_n)\ts$-modules by restriction.
Denote by $I$ the operator which acts as the composition
of \eqref{hint} with \eqref{lint} on the first tensor factor
the source $\Y(\g_n)\ts$-module in \eqref{intoper2}, and
which acts trivially on other $m-1$ tensor factors of the source module.

The intertwiner of $\Y(\g_n)\ts$-modules \eqref{bint}
can also be regarded as that of~the $\Y(\g_n)\ts$-modules
\begin{equation}
\label{cint}
\Phd^{\,n-\nu_1}_{\frac12-\mu_1-\rho_1}
\to\,
\Phd^{\,\nu_1-n}_{\frac12-\mu_1-\rho_1}
\end{equation}
where we use the notation introduced 
immediately before stating Lemma \ref{1.15}.
Denote by $J$ the operator which acts as the composition
of \eqref{lint} with \eqref{cint} on the first tensor factor
the target $\Y(\g_n)\ts$-module in \eqref{intoper2}, and
which acts trivially on other $m-1$ tensor factors of the target module.

By definition, $J(s_0\circ\mu)$ is an intertwining operator 
of $\Y(\g_n)\ts$-modules
\begin{align*}
\Ph^{\,n-\nu_1}_{\frac12-\mu_1-\rho_1}
\ot\ts
\Ph^{\,\nu_2}_{\mu_2+\rho_2+\frac12}
&\ot\ts\ldots\ts\ot\ts
\Ph^{\,\nu_m}_{\mu_m+\rho_m+\frac12}
\\
&\ts\downarrow
\\
\Ph^{\,\nu_1-n}_{\frac12-\mu_1-\rho_1}
\ot\ts
\Ph^{\,-\nu_2}_{\mu_2+\rho_2+\frac12}
&\ot\ts\ldots\ts\ot\ts
\Ph^{\,-\nu_m}_{\mu_m+\rho_m+\frac12}\ .
\end{align*}
Let us now replace the first tensor factors of
the above two $\Y(\g_n)\ts$-modules by
$$
\Phd^{\,n-\nu_1}_{\frac12-\mu_1-\rho_1}
\quad\text{and}\quad
\Phd^{\,\nu_1-n}_{\frac12-\mu_1-\rho_1}
$$
respectively. The operator $J(s_0\circ\mu)$
also intertwines the resulting two tensor products
as  $\Y(\g_n)\ts$-modules. Take $J(s_0\circ\mu)$ in
its latter capacity. 
Then arguing like in the end of the proof of
Lemma \ref{2.55}, that is either performing a direct calculation,
or using the irreducibility of the source and target
$\Y(\g_n)\ts$-modules in \eqref{intoper2} for any generic weight
$\mu$ of $\f_m\ts$, we obtain the relation
$$
J\,J(\mu)=J(s_0\circ\mu)\,I 
$$
for any dominant weights $\la+\rho$ and $s_0(\la+\mu)$ of $\f_m\ts$. 
This proves Lemma \ref{3.52}, because both $I$ and $J$
are invertible and intertwine $\Y(\g_n)\ts$-modules.
\qed
\end{proof}

For any $\la\in\h_m^*$ let $H_\la$ be
the subgroup of $H_m$ consisting of all elements $w$ such that 
$w\circ\la=\la\ts$. Let $\mathcal{O}$ be an orbit of the shifted 
action of the subgroup $H_\la\subset H_m$ on $\h_m^*\ts$.
If $\nu_1\lcd\nu_m\in\{\tts1\lcd n-1\}$ for at least one
weight $\mu\in\mathcal{O}$, then every $\mu\in\mathcal{O}$
satisfies the same condition. 
Suppose this is the case for $\mathcal{O}$.
If $\la+\rho$ is dominant,
then there is at least one weight $\mu\in\mathcal{O}$
such that the pair $(\la\com\mu)$ is good.
Theorem \ref{1.2} generalizes due to the following proposition.

\begin{proposition*}
\label{3.5}
If $\la+\rho$ is dominant, then for all 
$\mu\in\mathcal{O}$
the images of the corresponding operators\/ $J(\mu)$
are similar to each other as\/ $\Y(\g_n)\ts$-modules.
\end{proposition*}

\begin{proof}
Take any $w\in H_\la\ts$ and any reduced decomposition 
$w=s_{c_l}\ldots s_{c_1}\ts$. Here $c_1\lcd c_l\ge0\ts$.
It can be derived from \cite[Corollary VI.1.2]{B} that
the weight $s_{c_k}\ldots s_{c_1}(\la+\rho)$ of $\f_m$ is
dominant for each $k=1\lcd l\ts$.
Applying Lemmas \ref{3.51} and \ref{3.52} now
completes the proof of Proposition \ref{3.5}.
\qed
\end{proof}


Thus all assertions of Theorem \ref{1.2} will remain valid
if we replace the good pair there by any pair
$(\la\com\mu)$ such that the weight $\la+\rho$ of $\f_m$ is dominant.   
However, we still assume that $\nu_1\lcd\nu_m\in\{\tts1\lcd n-1\}$
for the latter pair.


\section*{\normalsize\bf Acknowledgements}

We are grateful to Ivan Cherednik for 
helpful discussions. The first author was supported by 
the RFBR grant 11-01-00962, by the joint RFBR-Royal Society 
grant 11-01-92612, and by the Federal Agency
for Science and Innovations of the Russian Federation
under the contract  14.740.11.0347.


\renewcommand\refname{{\normalsize\bf References}}



\end{document}